\documentclass{article}

\author{Dixy Msapato}
 \usepackage{amsmath,enumerate,amsthm,amsfonts,amssymb,latexsym,bbm}
\usepackage[super]{nth}
\usepackage[a4paper]{geometry}
\usepackage[hidelinks]{hyperref}
\usepackage[UKenglish]{babel}
\usepackage{float}
\usepackage{dirtytalk}
\usepackage{tikz}
\usepackage{tikz-qtree}
\usepackage{mathtools}
\usepackage{graphicx}
\usepackage[mathscr]{euscript} 

\newdimen\R
\newdimen\T
\newcommand\blfootnote[1]{%
  \begingroup
  \renewcommand\thefootnote{}\footnote{#1}%
  \addtocounter{footnote}{-1}%
  \endgroup
}

\theoremstyle{definition}
\newtheorem{theorem}{Theorem}[section]
\newtheorem{definition}[theorem]{Definition}
\newtheorem{proposition}[theorem]{Proposition}
\newtheorem{remark}[theorem]{Remark}
\newtheorem{corollary}[theorem]{Corollary}

\newtheorem{lemma}[theorem]{Lemma}

\newtheorem{example}[theorem]{Example}

\title{Modular Fuss-Catalan Numbers}
\date{}
\begin{document}
\maketitle
\begin{abstract}
The modular Catalan numbers $C_{k,n}$, introduced by Hein and Huang in 2016 count equivalence classes of parenthesizations of $x_0 * x_1 * \dots *x_n$ where $*$ is a binary $k$-associative operation and $k$ is a positive integer. The classical notion of associativity is just 1-associativity, in which case $C_{1,n} = 1$ and the size of the unique class is given by the Catalan number $C_n$. In this paper we introduce modular Fuss-Catalan numbers $C_{k,n}^{m}$ which count equivalence classes of parenthesizations of $x_0 * x_1 * \dots *x_n$ where $*$ is an $m$-ary $k$-associative operation for $m \geq 2$. Our main results are a closed formula for $C_{k,n}^{m}$ and a characterisation of $k$-associativity. 
\end{abstract}



\thanks{}
\blfootnote{Address: School of Mathematics, University of Leeds, Leeds, LS2 9JT, United Kingdom. \\
	Contact: mmdmm@leeds.ac.uk \\
	\thanks{This research was supported by an EPSRC Doctoral Training Partnership (reference EP/R513258/1) through the University of Leeds. The author also wishes to thank their supervisor, Robert Marsh for their constant support.}}
\section{Introduction}
Let $X$ be a non-empty set with a binary operation $\star:X^2 \rightarrow X$ and let $n \geq 1$ be a positive integer. If $\star$ is associative, the general associativity law states that the expression $x_1 \star x_2 \star \dots \star x_n$ is unambiguous  for all $x_1,x_2, \dots, x_n \in X$, i.e. all possible parenthesizations of the expression result in the same evaluation. We say a non-associative operation $\star$ is \textit{left-associative} if the order of operation is understood to be from left to right, in which case we write $x_1 \star \dots \star x_n$ to mean $((\dots ((x_1 \star x_2) \star x_3) \dots \star x_{n-1})\star x_n).$ Let $k \geq 1$ be a positive integer. There is a notion of $k$-associativity for binary operations which generalises the notion of associativity. A left-associative binary operation $\star$ is \textit{$k$-associative} if $$(x_1 \star \dots \star x_{k+1})\star x_{k+2} = x_1 \star (x_2 \star \dots \star x_{k+1} \star x_{k+2}) \text{ for all } x_{1}, x_{2}, \dots, x_{k+2} \in X.$$
By setting $k=1$, we recover the classical notion of associativity for binary operations. In the case where $k>1$, the general associativity law no longer holds, that is in general the evaluation of the expression $x_1 \star x_2 \star \dots \star x_n$ depends on its parenthesization. The $k$-associative binary operations are studied in \cite{HeinHuang}.

Fix a positive integer $m \geq 2$. An \textit{$m$-ary operation} on $X$ is a map $* : X^{m} \rightarrow X$. Another way to generalise associativity of binary operations is to consider associative $m$-ary operations. We say that the $m$-ary operation $*$ is associative if for $1 \leq j \leq m-1,$
 \begin{multline}\label{eqn:m-assoc}
x_1 * \dots *x_{j-1}*(x_j*x_{j+1}*\dots*x_{j+(m-1)})*x_{j+(m-1)+1}*x_{j +(m-1)+2}*\dots*x_{m+(m-1)} \\
 = x_1*\dots*x_{j-1}*x_{j}*(x_{j+1}*\dots*x_{j+(m-1)}*x_{j+(m-1)+1})*x_{j+(m-1)+2}*\dots*x_{m+(m-1)}.
\end{multline}
As in the case for associative binary operations, we have a general associativity law stating that the expression $x_1 * x_2 * \dots * x_n$ is independent of parenthesization (for example, see \cite[ Theorem 2.1]{Andres}). We note that $n$ is not arbitrary in this case, but is of the form $n = m+g(m-1)$ for some integer $g \geq 1$. Associative $m$-ary operations are important for the study of $m$-semigroups and polyadic groups. These are generalisations of semigroups and groups where we consider associative $m$-ary operations instead of associative binary operations. The $m$-semigroups were introduced in \cite{Dornte} and polyadic groups were introduced extensively in \cite{Post} and \cite{UniversalAlgebras}.

In this paper we will study $m$-ary $k$-associative operations, which are a further generalisation of associative binary operations that combines the two above generalisations. We say $*$ is left-associative if the order of operation is from left to right, meaning for an integer $g \geq 1$, we write $x_1 * x_2 * \dots * x_{m+g(m-1)}$ to mean $$((\dots((x_1 *  \dots * x_m)*x_{m+1} * \dots * x_{m+(m-1)})\dots * x_{m+(g-1)(m-1)})x_{m+(g-1)(m-1)+1} \dots x_{m+g(m-1))}.$$ A left-associative $m$-ary operation $*$ is said to be $k$-associative if for $1 \leq j \leq m-1,$ the following equality holds:
\begin{multline}\label{eqn:1}
x_1* \dots *x_{j-1}*(x_j*x_{j+1}*\dots*x_{j+k(m-1)})*x_{j+k(m-1)+1}*x_{j+k(m-1)+2}*\dots*x_{m+k(m-1)} \\
 = x_1*\dots*x_{j-1}*x_{j}*(x_{j+1}*\dots*x_{j+k(m-1)}*x_{j+k(m-1)+1})*x_{j+k(m-1)+2}*\dots*x_{m+k(m-1)}.
\end{multline}
We note that the terminology \say{$k$-associativity} is used by Wardlaw in \cite{wardlaw2001generalized} to mean associativity of $k$-ary operations. This is not to be confused with the notion of $k$-associativity we consider here, which is a generalisation of associativity for $m$-ary operations (and coincides with the $k$-associativity for binary operation introduced in \cite{HeinHuang}). Let $*$ be a $k$-associative $m$-ary operation. Let $g \geq 1$ be a positive integer. Then for $n =m+g(m-1)$ and $k >1$, the expression $x_1*\dots*x_n$ is ambiguous without a parenthesization to clarify the order of operation i.e. the general associativity law no longer holds. Let $p$ and $p^{\prime}$ be two parenthesizations of $x_1*\dots*x_n$.  If we may obtain $p^{\prime}$ from $p$ by finitely many left-hand-side to right-hand-side applications of the $k$-associative property (\ref{eqn:1}) to $p$, we write $p \preceq_{k} p^{\prime}$. This induces a partial order on the set of parenthesizations of $x_1 * x_2 * \dots * x_n$ called the \textit{$k$-associative order}. The connected components of the $k$-associative order are called \textit{$k$-components}. We say that two parenthesizations of $x_1 * x_2 * \dots * x_n$ are \textit{$k$-equivalent} if they lie in the same $k$-component. When $k=1$ and $m=2$, we recover the well-known \textit{Tamari lattice} (see for example \cite{geyer1994tamari}). In general, determining whether two parenthesizations are $k$-equivalent is a non-trivial problem. 

Let $A =( \mathbb{C} \langle u_1, \dots , u_n \rangle, +,  \cdot)$ be the free associative algebra over $\mathbb{C}$ in $n$ indeterminates $u_1,u_2, \dots, u_n$. Let $\omega$ be an element of order $k(m-1)$ in $A$. We define an $m$-ary operation $\circ : A^m \rightarrow A$ in the following way for $a_1, \dots,a_m $ in A: $$a_1 \circ \dots \circ a_m = \omega^{m-1} \cdot a_1 + \omega^{m-2} \cdot a_2 + \dots \omega \cdot a_{m-1} + a_m.$$ It is easy to show that this operation is $k$-associative by direct calculation. 

\begin{theorem}\label{mainThmIntro}
Let $*:X^{m} \rightarrow X$ be a $k$-associative $m$-ary operation on a set $X$ where $m\geq 2$ and $k \geq 1$ are integers. Let $p(x_1*\dots*x_n) \text{ and }p^{\prime}(x_1 * \dots *x_n)$ be two $m$-ary parenthesizations of the $m$-ary expression $x_1*\dots*x_n$, where $n = m+g(m-1)$ for some positive integer $g$. Then $p(x_1*\dots*x_n)$ is $k$-equivalent to $p^{\prime}(x_1 * \dots *x_n)$ if and only if  $$p(u_1 \circ \dots \circ u_n)= p^{\prime}(u_1 \circ \dots \circ u_n),$$ where the equality comes from evaluating the parenthesizations under $\circ$ in the algebra $A$. 
\end{theorem}
We define the \textit{$(k)$-modular Fuss-Catalan number} $C^{m}_{k,n}$ to be the number of $k$-equivalence classes of parenthesizations of $x_0*x_1*\dots*x_n$. By the general associativity law we have that $C^{m}_{1,n} = 1$, and the size of this class is given by the Fuss-Catalan number $C^m_n = \frac{1}{(m-1)n+1}\binom{mn}{n}$.
The following theorem is a generalisation of \cite[Proposition 2.10]{HeinHuang}.
\begin{theorem}\label{closedFormula}
The $(k)$-modular Fuss-Catalan number is given by the following closed formula, 
$$C^{m}_{k,n} = \sum_{\substack{ 1 \leq l \leq n \\ m-1 | l}} { \frac{l}{n} \sum_{\substack{ m_1 + \dots +m_k = n \\ m_2 + 2m_3 + \dots (k-1)m_k = \frac{n-l}{m-1}}}}{\binom{n}{m_1, \dots m_k}}.$$
\end{theorem}

This paper is organised as follows: in \S 2 we define right $k$-rotations, left $k$-rotations and $k$-equivalence for $m$-ary trees. In \S 3 we study $k$-equivalence by appealing to $m$-Dyck paths in order to prove Theorem \ref{thm:main}; which can be thought of as the $m$-ary tree version of Theorem \ref{mainThmIntro}. In \S 4, we prove our first main result Theorem \ref{mainThmIntro} using some results obtained in the previous section. Finally in \S 5 we derived the closed formula in Theorem \ref{closedFormula} using $m$-Dyck paths; a method adopted from \cite[\S 5]{HeinHuangv2}. 
\section{$m$-ary Trees}
In studying $k$-equivalence, it is often more convenient to do so by appealing to other sequences of combinatorial sets counted by the Fuss-Catalan numbers. In this section we will study $k$-equivalence via $m$-ary trees. In order to do this, we use a known bijection between parenthesizations of $m$-ary expressions and $m$-ary trees outlined in \cite[\S 0]{HiltonPedersen}. For the rest of this section, we fix integers $m \geq 2$, $g \geq 0$, $k \geq 1$ and $n = m+g(m-1)$.

\begin{definition}{\textbf{$m$-ary Tree }}\cite[Section 4, A14(b)]{Stanley}.
An $m$-ary tree is a rooted tree with the property that each node either has $0$ or $m$ linearly ordered children. A node with no children will be referred to as a \textit{leaf} and the unique node without a parent is called the  \textit{root}.
\end{definition}
\begin{remark}
The objects we are calling $m$-ary trees in this paper are commonly referred to as \textit{full $m$-ary trees} in the wider literature.
\end{remark} It will be our convention in this paper to draw $m$-ary trees with the root at the top and leaves at the bottom. We shall denote the set of $m$-ary trees with $n$ leaves by $\mathbf{B}_{n}^{m}$. We enumerate the leaves from left to right in a linear order starting from $1$ on the leftmost leaf, up to $n$ on the rightmost leaf. We will endow the $m$-ary trees with an additional edge labelling with labels drawn from the set $\{l_1, l_2, \dots, l_m\}$.  An edge will be given the label $l_i$ if it links a node with its $i^{\text{th}}$ child. See the figure below for an example. 

\begin{figure}[H]
\begin{center}
\begin{tikzpicture}

\draw[fill=black] (2,0) circle (1pt);
\draw[fill=black] (4,0) circle (1pt);
\draw[fill=black] (6,0) circle (1pt);
\draw[fill=black] (2,2) circle (1pt);
\draw[fill=black] (4,2) circle (1pt);
\draw[fill=black] (6,2) circle (1pt);
\draw[fill=black] (4,4) circle (1pt);
\draw[fill=black] (6,4) circle (1pt);
\draw[fill=black] (8,4) circle (1pt); 
\draw[fill=black] (6,6) circle (1pt); 

\draw[thin] (2,0) -- (4,2)node[midway, left] {$l_1$};
\draw[thin] (4,0) -- (4,2) node[midway, right] {$l_2$};
\draw[thin] (6,0) -- (4,2) node[midway, right] {$l_3$};
\draw[thin] (2,2) -- (4,4) node[midway, left] {$l_1$};
\draw[thin] (4,2) -- (4,4)node[midway, right] {$l_2$};
\draw[thin] (6,2) -- (4,4)node[midway, right] {$l_3$};
\draw[thin] (4,4) -- (6,6)node[midway, left] {$l_1$};
\draw[thin] (6,4) -- (6,6)node[midway, right] {$l_2$};
\draw[thin] (8,4) -- (6,6)node[midway, right] {$l_3$};

\node at (2,-0.4) {$2$};
\node at (4,-0.4) {$3$};
\node at (6,-0.4) {$4$};
\node at (2,1.6) {$1$};
\node at (6,3.6) {$6$};
\node at (6,1.6) {$5$};
\node at (8,3.6) {$7$};

\end{tikzpicture}
\end{center}
\caption{A labelled 3-ary tree.}\label{fig:3tree}
\end{figure}
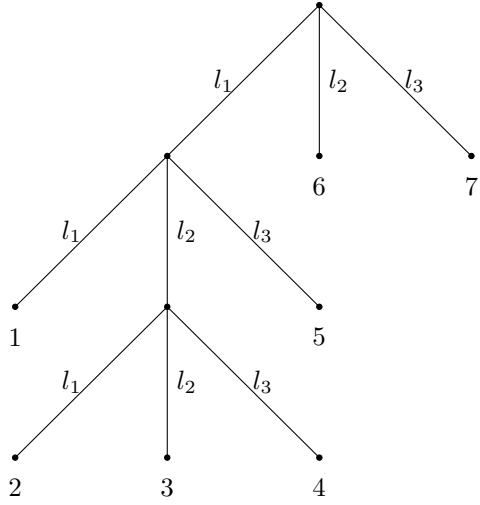

\begin{definition}{\textbf{Meet.}}
Let $t_1,t_2,\dots,t_m$ be $m$-ary trees. We define the \textit{meet} of $t_1,t_2,\dots,t_m$ to be the $m$-ary tree $t_1 \wedge t_2 \wedge\dots\wedge t_m$, which has the tree $t_i$ as the subtree rooted at the $i^{\text{th}}$ child of the root for $1 \leq i \leq n$. 
\end{definition}
We regard the meet as a left-associative $m$-ary operation, meaning that the order of operation is from left to right. That is to say, we will write $t_1 \wedge t_2 \wedge \dots \wedge t_{m+g(m-1)}$ to mean $$((\dots((t_1 \wedge  \dots \wedge t_m)t_{m+1} \wedge \dots \wedge t_{m+(m-1)})\dots \wedge t_{m+(g-1)(m-1)})t_{m+(g-1)(m-1)+1} \dots t_{m+g(m-1))}.$$
The following bijection is well-known, see for example \cite[\S 0]{HiltonPedersen}.
\begin{proposition}\label{HPBijection}\cite[\S 0]{HiltonPedersen}
Let $X$ be a non-empty set and let $*:X^{m} \rightarrow X$ be an $m$-ary operation. Take $x_1, x_2, \dots, x_n$ in $X$. Then there is a bijection between the set of $m$-ary trees on $n$ leaves and the set of $m$-ary parenthesizations of the expression $x_1 * x_2 * \dots * x_n$ which is defined in the following way. Let $t$ be an $m$-ary tree with $n$ leaves such that the $i^{\text{th}}$ leaf of $t$ is labelled $\varepsilon_i$. Consider the tree $t$ expressed as a bracketed meet of its leaves $\varepsilon_i$, where the $\varepsilon_i$ are thought of as trees consisting of just a root. The bijection maps $t$ to the parenthesization obtained by replacing $\wedge$ with $*$ and replacing $\varepsilon_i$ with $x_i$. The inverse map from the set of $m$-ary parenthesizations of the expression $x_1 * x_2 * \dots * x_n$ to $m$-ary trees with $n$ leaves acts in the naturally opposite way.
\end{proposition}

\begin{example}
Let $t$ be the tree in Figure \ref{fig:3tree}. Thinking of the leaves of $t$ as 3-ary trees consisting of just a root, assign to each leaf $i$ the label $\varepsilon_i$. We can write $t$  as a meet of the leaves $\varepsilon_{i}$ so $t=(\varepsilon_1 \wedge (\varepsilon_2 \wedge \varepsilon_3 \wedge \varepsilon_4) \wedge \varepsilon_5) \wedge \varepsilon_6 \wedge \varepsilon_7$. Under the bijection in Proposition $\ref{HPBijection}$ the tree $t$ is assigned to the parenthesization of the $x_i$ given by $(x_1*(x_2*x_3*x_4)*x_5)*x_6*x_7$.
\end{example}

\begin{definition}{\textbf{Right $k$-rotation.}}
Let $k \geq 1$ be a positive integer. Let $t_1, t_2, \dots, t_{(m-1)+k(m-1)}$ be $m$-ary trees.  Let $1 \leq j \leq m-1$. Suppose that $t \in \mathbf{B}_{n}^{m}$ has a subtree,
 $$s = t_1 \wedge t_2 \wedge\dots\wedge t_{j-1} \wedge(t_j \wedge t_{j+1} \wedge\dots\wedge t_{j+k(m-1)}) \wedge t_{j+k(m-1)+1} \wedge t_{j+k(m-1)+2} \wedge \dots \wedge t_{(m-1)+k(m-1)}$$
  rooted at some node $v$ in $t$. The \textit{right $k$-rotation} of $t$ at $v$ is the operation of replacing $s$ with the subtree $$s^{\prime} = t_1 \wedge t_2 \wedge\dots\wedge t_{j-1} \wedge t_{j} \wedge (t_{j+1} \wedge \dots\wedge t_{j+k(m-1)} \wedge t_{j+k(m-1)+1}) \wedge t_{j+k(m-1)+2} \wedge\dots\wedge t_{(m-1)+k(m-1)}.$$   
\end{definition}
\begin{remark}\label{rot&assoc}
It should be clear that under the bijection in Proposition \ref{HPBijection}, right $k$-rotation of $m$-ary trees corresponds to a left-hand-side to right-hand-side application of the $k$-associative rule in (\ref{eqn:1}).
 We can also define a left $k$-rotation dually by switching the roles of $s$ and $s^{\prime}$ in the above definition. In this case, a left $k$-rotation corresponds to a right-hand-side to left-hand-side application of the $k$-associative rule in (\ref{eqn:1}). \end{remark}

\begin{definition}Let $t$ and $t^{\prime}$ be $m$-ary trees with $n$ leaves. If we can obtain $t^{\prime}$ from $t$ by applying finitely many right $k$-rotations to $t$, we write $t \preceq_{k} t^{\prime}$. The induced partial order on $\mathbf{B}_{n}^{m}$ is called the $k$-\textit{associative} order. The connected components (i.e connected components of the Hasse diagram) of $\mathbf{B}_{n}^{m}$ under the $k$-associative order are called $k$-\textit{components}. We then say two $m$-ary trees with $n$ leaves are $k$-\textit{equivalent} if they belong to the same $k$-component of $\mathbf{B}_{n}^{m}$.
\end{definition}

\begin{example}{\label{rotex}}
The example below shows a right 2-rotation on a 3-ary tree. We apply the right 2-rotation at $v$.\\
\begin{figure}[H]
\begin{center}
\begin{tikzpicture}
 
\draw[fill=black] (0,0) circle (1pt);
\draw[fill=black] (1,0) circle (1pt);
\draw[fill=black] (2,0) circle (1pt);
\draw[fill=black] (1,1) circle (1pt);
\draw[fill=black] (2,1) circle (1pt);
\draw[fill=black] (3,1) circle (1pt);
\draw[fill=black] (2,2) circle (1pt);
\draw[fill=black] (3,2) circle (1pt);
\draw[fill=black] (4,2) circle (1pt); 
\draw[fill=black] (3,3) circle (1pt);

\draw[thin] (0,0) -- (1,1);
\draw[thin] (1,0) -- (1,1);
\draw[thin] (2,0) -- (1,1);
\draw[thin] (1,1) -- (2,2);
\draw[thin] (2,1) -- (2,2);
\draw[thin] (3,1) -- (2,2);
\draw[thin] (2,2) -- (3,3);
\draw[thin] (3,2) -- (3,3);
\draw[thin] (4,2) -- (3,3);

\node at (0,-0.2) {$t_1$};
\node at (1,-0.2) {$t_2$};
\node at (2,-0.2) {$t_3$};
\node at (2,0.8) {$t_4$};
\node at (3,0.8) {$t_5$};
\node at (0.7,1) {};
\node at (3,1.8) {$t_6$};
\node at (4,1.8) {$t_7$};
\node at (1.7,2) {};
\node at (3,3.2) {$v$};

\draw[fill=black] (7,0) circle (1pt);
\draw[fill=black] (8,0) circle (1pt);
\draw[fill=black] (9,0) circle (1pt);
\draw[fill=black] (8,1) circle (1pt);
\draw[fill=black] (10,1) circle (1pt);
\draw[fill=black] (9,1) circle (1pt);
\draw[fill=black] (10,1) circle (1pt);
\draw[fill=black] (9,2) circle (1pt);
\draw[fill=black] (8,2) circle (1pt);
\draw[fill=black] (10,2) circle (1pt); 
\draw[fill=black] (9,3) circle (1pt); 

\draw[thin] (7,0) -- (8,1);
\draw[thin] (8,0) -- (8,1);
\draw[thin] (9,0) -- (8,1);
\draw[thin] (8,1) -- (9,2);
\draw[thin] (9,1) -- (9,2);
\draw[thin] (10,1) -- (9,2);
\draw[thin] (8,2) -- (9,3);
\draw[thin] (9,2) -- (9,3);
\draw[thin] (10,2) -- (9,3);

\node at (7,-0.2) {$t_2$};
\node at (8,-0.2) {$t_3$};
\node at (9,-0.2) {$t_4$};
\node at (8,1.8) {$t_1$};
\node at (9,0.8) {$t_5$};
\node at (7.7,1) {};
\node at (10,0.8) {$t_6$};
\node at (10,1.8) {$t_7$};
\node at (9.2,2.1) {};
\node at (9,3.2) {$v$};

\end{tikzpicture}
\end{center}
\caption{The tree on the right is a result of a right 2-rotation of the tree on the  left at $v$.}
\end{figure}
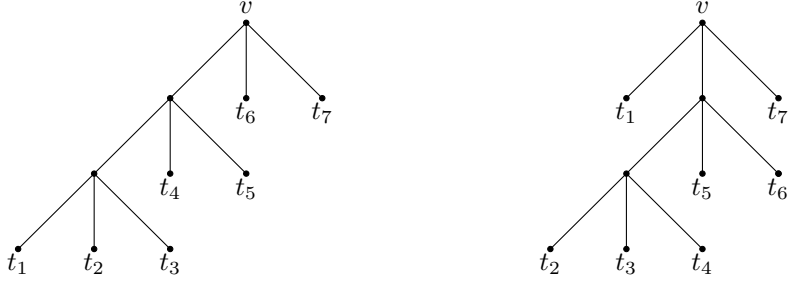

The subtree rooted at $v$ is $s = (t_1 \wedge t_2 \wedge t_3 \wedge t_4 \wedge t_5 ) \wedge t_6 \wedge t_7 =   ((t_1 \wedge t_2 \wedge t_3) \wedge t_4 \wedge t_5 ) \wedge t_6 \wedge t_7$. The subtree $s$ is then replaced  by the subtree $s^{\prime} =t_1 \wedge (t_2 \wedge t_3 \wedge t_4 \wedge t_5 \wedge t_6) \wedge t_7 = t_1 \wedge ((t_2 \wedge t_3 \wedge t_4) \wedge t_5 \wedge t_6) \wedge t_7$ at $v$. 
\end{example}
The following proposition is a generalisation of \cite[Proposition 2.5]{HeinHuang}. 
\begin{proposition}
Let $t$ be an $m$-ary tree such that we can perform a right $k$-rotation of $t$ at some node $v$. Suppose that $k=pk^{\prime}$ for some positive integers $p$ and $k^{\prime}$. Then the right $k$-rotation at $v$ may be decomposed into a sequence of $p$ right $k^{\prime}$-rotations of $t$. The same holds for left $k$-rotations
\begin{proof}
We argue by induction on $p$. The case for $p=1$ is trivial. So assume for induction that the statement is true for some $p \geq 1$. Suppose that $k=(p+1)k^{\prime}$ for some positive integer $k^{\prime}$.

Suppose we have a tree $t$ which we can right $k$-rotate at some node $v$. Denote by $s$ the subtree of $t$ rooted at $v$. For some $1 \leq j \leq m-1$, 

$s = t_1 \wedge t_2 \wedge \dots \wedge t_{j-1} \wedge (t_j \wedge t_{j+1} \wedge \dots \wedge t_{j+pk^{\prime}(m-1)} \wedge t_{j+pk^{\prime}(m-1) + 1} \wedge \dots \wedge t_{j+k(m-1)})\wedge t_{j+k(m-1)+1} \wedge \dots \wedge t_{(m-1)+k(m-1)}.$

The right $k$-rotation replaces the subtree $s$ with the subtree $s^{\prime}$ where \\
$s^{\prime} =   t_1 \wedge t_2 \wedge \dots \wedge t_{j-1} \wedge t_j \wedge (t_{j+1} \wedge \dots \wedge t_{j+m-1} \wedge t_{j+m} \wedge \dots \wedge t_{j+pk^{\prime}(m-1)} \wedge t_{j+pk^{\prime}(m-1) + 1} \wedge t_{j + pk^{\prime}(m-1)+2} \wedge\dots\wedge t_{j+k(m-1)} \wedge t_{j+k(m-1)+1}) \wedge\dots\wedge t_{(m-1)+k(m-1)}.$

We will show that the result of this right $k$-rotation can also be obtained by performing $(p+1)$ right $k^{\prime}$-rotations. 

Let $r$ be the following subtree of $s$, which is rooted at the $j^{\text{th}}$ child of the root of $s$, $$r = (t_j \wedge t_{j+1} \wedge \dots \wedge t_{j+pk^{\prime}(m-1)} \wedge t_{j+pk^{\prime}(m-1) + 1} \wedge \dots \wedge t_{j+k(m-1)}).$$
We can write $$r=((t_j \wedge t_{j+1} \wedge \dots \wedge t_{j+m-1} \wedge t_{j+m} \wedge \dots \wedge t_{j+pk^{\prime}(m-1)}) \wedge t_{j+pk^{\prime}(m-1) + 1} \wedge \dots \wedge t_{j+k(m-1)})$$
since the meet operation is left-associative.
Performing a right $pk^{\prime}$-rotation of $t$ at the $j^{\text{th}}$ child of the root of $s$, we replace $r$ with $$r^{\prime} = (t_j \wedge (t_{j+1} \wedge \dots \wedge t_{j+m-1} \wedge t_{j+m} \wedge \dots \wedge t_{j+pk^{\prime}(m-1)} \wedge t_{j+pk^{\prime}(m-1) + 1}) \wedge \dots \wedge t_{j+k(m-1)}).$$
By the inductive hypothesis, this right $pk^{\prime}$-rotation is the result of $p$ right $k^{\prime}$-rotations. 
  
Set $$u := (t_{j+1} \wedge \dots \wedge t_{j+m-1} \wedge t_{j+m} \wedge \dots \wedge t_{j+pk^{\prime}(m-1)} \wedge t_{j+pk^{\prime}(m-1) + 1}).$$ Then $$r^{\prime} = (t_j \wedge u \wedge \dots \wedge t_{j+k(m-1)}).$$
Thus the above right $pk^{\prime}$-rotation of $t$ at the $j^{\text{th}}$ child of the root of $s$, replaces $s$ with $q$ at the node $v$ in $t$ where,  $$q = t_1 \wedge t_2 \wedge \dots \wedge t_{j-1} \wedge (t_j \wedge u \wedge t_{j + pk^{\prime}(m-1)+2} \wedge\dots\wedge t_{j+k(m-1)}) \wedge t_{j+k(m-1)+1} \wedge\dots\wedge t_{(m-1)+k(m-1)}.$$ 

We then perform a right $k^{\prime}$-rotation of $t$ at $v$. This replaces $q$ with the subtree, $$q^{\prime} =  t_1 \wedge t_2 \wedge \dots \wedge t_{j-1} \wedge t_j \wedge (u \wedge t_{j + pk^{\prime}(m-1)+2} \wedge\dots\wedge t_{j+k(m-1)} \wedge t_{j+k(m-1)+1}) \wedge\dots\wedge t_{(m-1)+k(m-1)}.$$ It is easy to see that $s^{\prime} = q^{\prime}$, therefore the result of performing the right $k=(p+1)k^{\prime}$-rotation at $v$ is precisely the result of performing $(p+1)$ right-$k^{\prime}$ rotations. The proof for left $k$-rotations is similar. 
\end{proof}
\end{proposition}

\begin{definition}{\textbf{Path.}}
Let $t$ be an $m$-ary tree and $n$ a positive integer. A \textit{path} $p$ in $t$ of length $n$ from a node $v$ to a node $w$ is a sequence $p = (v_1, v_2, \dots , v_n)$ of nodes such that $v_1 = v, v_n =w$ and $(v_i,v_{i+1})$ is an edge in $t$ for $1 \leq i \leq n-1$.
\end{definition}

\begin{definition}{\textbf{Depth.}}
Let $t$ be an $m$-ary tree with $n$ leaves and edges labelled by labels drawn from the set $\{l_1, l_2, \dots, l_m\}$. For $i=1,2, \dots , m$ and $j=1,2, \dots, n$, let $\delta_{j}^{l_i}(t)$ be the number of edges labelled $l_i$ in the unique path from the root to the $j^{\text{th}}$ leaf. Let $\delta^{l_i}(t) = (\delta_{1}^{l_i}(t), \dots , \delta_{n}^{l_i}(t))$ and set $\delta(t) = (\delta^{l_1}(t), \dots , \delta^{l_n}(t))$. We call $\delta(t)$ the depth of $t$.  
\end{definition}
\begin{example}
Let $t$ be the tree in Figure \ref{fig:3tree}. The depth of tree $t$ is given by $$\delta(t) = ((2,2,1,1,0,0,0),(0,1,2,1,0,1,0),(0,0,0,1,1,0,1)).$$
\end{example}
The following lemmas are easy to verify. 
\begin{lemma}\label{rem1}
Let $t$ be an $m$-ary tree with $n$ leaves and depth $( \delta^{l_1}(t), \delta^{l_2}(t),\dots,\delta^{l_{m}}(t) )$. We have that $\delta_{n}^{l_m}(t)  \neq 0$ and $\delta_{n}^{l_i}(t) = 0$ for $i \neq m$. Dually, $\delta_{1}^{l_1}(t)  \neq 0$ and $\delta_{1}^{l_i}(t) = 0$ for $i \neq 1$. 
\end{lemma}
This is because the unique path from the root to the $n^{\text{th}}$ leaf involves  choosing the $m^{\text{th}}$ child at each stage. Similarly for the dual statement. 

\begin{lemma}\label{rem2}
Let $t$ be an $m$-ary tree with $n$ leaves and depth $( \delta^{l_1}(t), \delta^{l_2}(t),\dots,\delta^{l_{m}}(t) )$. We have that $\delta_{n-1}^{l_{m-1}}(t) = 1$. For $1 \leq i \leq m-2$, we have that $\delta_{n-1}^{l_i}(t)=0$.
\end{lemma}
This is because the unique path from the root to the $(n-1)^{\text{th}}$ leaf involves choosing the $m^{\text{th}}$ child at every stage but one, in which case, we choose the $(m-1)^{\text{th}}$ child.
 
We shall prove the following result in the next section.
\begin{theorem}\label{thm:main}
Let $t,t^{\prime}$ be a pair of $m$-ary trees with $n$ leaves with depths $(\delta^{l_1}(t), \delta^{l_2}(t),\dots,\delta^{l_{m}}(t) )$ and $(\delta^{l_1}(t^{\prime}), \delta^{l_2}(t^{\prime}),\dots,\delta^{l_{m}}(t^{\prime}))$ respectively. We have that $t$ and $ t^{\prime}$ are $k$-equivalent if and only if \[  \sum\limits_{i=1}^{m-1} (m-i)\delta^{l_{i}}(t) \equiv \sum\limits_{i=1}^{m-1} (m-i)\delta^{l_{i}}(t^{\prime}) \text{ mod } k(m-1)  \] where the addition on the $n$-tuples is componentwise. 
\end{theorem}
The strength of the theorem is that it allows us to determine the $k$-equivalence of $m$-ary trees from simply reading their depths. 
The case $m=2$ is known, see \cite[Proposition 2.11]{HeinHuang}). We shall prove the case for general $m \geq 2$. To do this, we appeal to another sequence of combinatorial sets counted by the Fuss-Catalan numbers, the $m$-Dyck paths. 

\section{m-Dyck Paths}
In this section we prove Theorem \ref{thm:main}. In order to do so, we appeal to a generalisation of Dyck paths known as $m$-Dyck paths to further study $k$-equivalence. We prove the theorem by first proving an $m$-Dyck path version of it. For the rest of this section, we fix integers $m \geq 2$, $g \geq 0$, $k \geq 1$ and $n = m+g(m-1)$.

\begin{definition}{\textbf{$m$-Dyck Path.}}
An \textit{$m$-Dyck path} is a lattice path in $\mathbb{Z}^2$ starting at $(0,0)$ consisting of up-steps $(m,m)$ and down-steps $(1,-1)$, which remains above the $x$-axis and ends on the $x$-axis. The length of a Dyck path is defined to be the number of down-steps it has.
\end{definition}

\begin{definition}{\textbf{Translated $m$-Dyck Path.}}
Let $a,b$ be non-negative integers both not equal to 0. A \textit{translated $m$-Dyck path} is a lattice path in $\mathbb{Z}^2$ starting at the point $(a,b)$ consisting of up-steps $(m,m)$ and down-steps $(1,-1)$, which remains above the line $y=b$ and ends on the line $y=b$.
\end{definition}

We denote the set of $m$-Dyck paths of length $n$ by $\mathbf{D}_{n}^{m}$. Where it is convenient, we refer to these paths as Dyck paths instead of $m$-Dyck paths. When referring to a translated $m$-Dyck path that is a sub path of a larger $m$-Dyck path, we will call it a \textit{sub $m$-Dyck path} or just \textit{sub-Dyck path}.  The following lemma is straight forward, so we state it without proof.

\begin{lemma}\label{form}
Let $D$ be an $m$-Dyck path of length $n$. Then we can write $$D = N^{d_1}SN^{d_2}S\dots SN^{d_n}S,$$ where $N$ denotes the up-step $(1,1)$ and $S$ denotes the down-step (1,-1). Note that when $m \neq 1$ the up-steps (1,1) are not steps on the path $D$ since by definition up-steps of $D$ are of the form $(m,m)$. Here $N^{d_i}$ is taken to mean a sequence of $d_i$ consecutive up-steps $N$. The $d_i$ are non-negative integer multiples of $m$ such that $d_1 + \dots + d_n = n$. Moreover the  $n$-tuple $d(D)=(d_1,d_2, \dots, d_n)$ is unique to each $m$-Dyck path $D$. 
 \end{lemma}

\begin{figure}[H]
\begin{center}

\begin{tikzpicture}
\draw[thin] (0,0)--(3,3)--(4,2)--(7,5)--(12,0);
\draw[fill=black] (0,0) circle (1pt);
\draw[fill=black] (3,3) circle (1pt);
\draw[fill=black] (4,2) circle (1pt);
\draw[fill=black] (7,5) circle (1pt);
\draw[fill=black] (8,4) circle (1pt);
\draw[fill=black] (9,3) circle (1pt);
\draw[fill=black] (10,2) circle (1pt);
\draw[fill=black] (11,1) circle (1pt);
\draw[fill=black] (12,0) circle (1pt);
\end{tikzpicture} 
\end{center}

\caption{a 3-Dyck path with 2 (3,3) up-steps. $D=N^{3}SN^{3}SN^{0}SN^{0}SN^{0}SN^{0}S$.}\label{fig:dyck}
\end{figure}
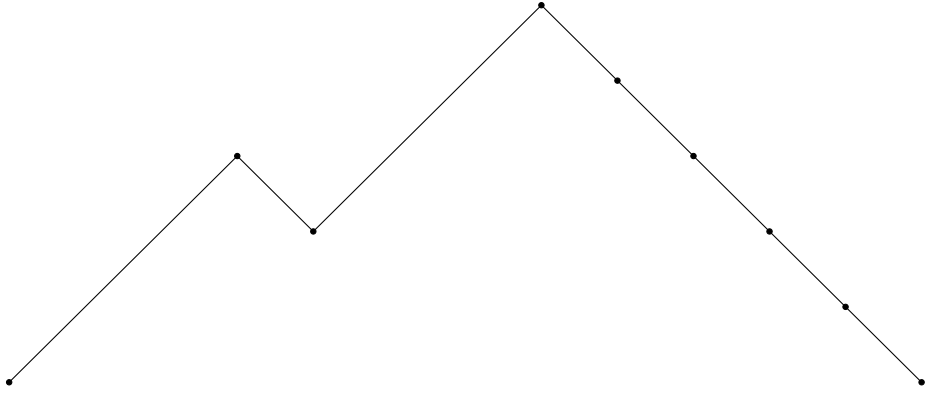
Let $D=  N^{d_1}SN^{d_2}S\dots SN^{d_n}S$ be an $m$-Dyck path of length $n$. When expressing $D$ in this way, if $d_i = 0$, we will omit $N^{d_i}$ from the expression. In this form we will also write $S^{l}$ to mean a sequence of $l$ consecutive $S$ steps. For example,  $D=N^{3}SN^{3}SN^{0}SN^{0}SN^{0}SN^{0}S = N^{3}SN^{3}S^{5}$.

We can express any $m$-ary tree as the meet of the $m$-ary sub-trees rooted at the children of the root. Therefore for $t$ an $m$-ary tree with $n$ leaves, we write $$t=t_1 \wedge t_2 \wedge \dots \wedge t_m,$$ 
where for $1 \leq i \leq m$ each $t_i$ is an $m$-ary tree with $n_i$ leaves and $n_1 + \dots + n_m = n$. 

Let $\varepsilon$ be the element of the singleton set $\mathbf{B}_{0}^{m}$, so $\varepsilon$ is the $m$-ary tree which consists of just a root. We construct a map  $\sigma_{m} \colon \mathbf{B}^{m} \to \mathbf{D}^{m-1}$ from the set of $m$-ary trees to the set of $(m-1)$-Dyck paths. We define $\sigma_{m}$ inductively in the following way, 
$$
\sigma_{m}(t) =
\begin{cases}
N^{0}S^{0} & \text{ if } t=\varepsilon ;\\
 N^{m-1} \sigma_{m}(t_1) S \sigma_{m}(t_2) S \dots S \sigma_{m}(t_m) & \text{otherwise}.\\
\end{cases}
$$
This construction generalises a well known map between binary trees (2-ary trees) and Dyck paths (1-Dyck paths); see for example \cite[Page 58, Tamari Lattice, Paragraph 2]{Bernadi}. 

\begin{example}{\label{Tree2Dyck}}
Consider the following 3-ary tree $t=\varepsilon \wedge \varepsilon \wedge ( \varepsilon \wedge \varepsilon \wedge \varepsilon)$. We calculate $\sigma_{3}(t),$
$$\sigma_{3}(t)=N^{2}\sigma_{3}(\varepsilon)S\sigma_{3}(\varepsilon)S\sigma_{3}(\varepsilon \wedge \varepsilon \wedge \varepsilon)$$
$$\hspace{1.95cm}= N^{2}N^{0}S^{0}SN^{0}S^{0}SN^2\sigma(\varepsilon)S\sigma(\varepsilon)S\sigma(\varepsilon)$$
$$\hspace{-1.9cm}= N^{2}S^{2}N^{2}S^{2}.
$$
See Figure \ref{fig:sigma}.

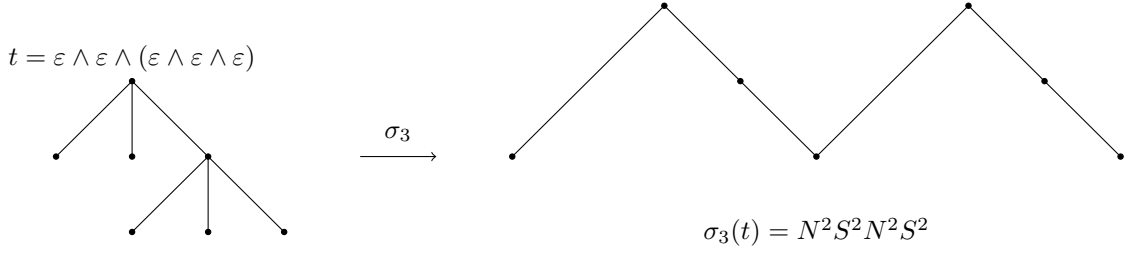
\begin{figure}[H]
\begin{center}
\begin{tikzpicture}
\draw[fill=black] (0,0) circle (1pt);
\draw[fill=black] (-1,-1) circle (1pt);
\draw[fill=black] (0,-1) circle (1pt);
\draw[fill=black] (1,-1) circle (1pt);
\draw[fill=black] (0,-2) circle (1pt);
\draw[fill=black] (1,-2) circle (1pt);
\draw[fill=black] (2,-2) circle (1pt);

\draw[thin] (0,0) -- (-1,-1);
\draw[thin] (0,0) -- (0,-1);
\draw[thin] (0,0) -- (1,-1);
\draw[thin] (1,-1) -- (0,-2);
\draw[thin] (1,-1) -- (1,-2);
\draw[thin] (1,-1) -- (2,-2);

\draw[thin,->] (3,-1)--(4,-1);
\node at (3.5,-0.7) {$\sigma_{3}$};
\node at (0,0.3) {$t = \varepsilon \wedge \varepsilon \wedge ( \varepsilon \wedge \varepsilon \wedge \varepsilon)$};
\node at (-1,-1.3){};
\node at (0,-1.3){};
\node at (2.3,-1){};
\node at (0,-2.3){};
\node at (1,-2.3){};
\node at (2,-2.3){};

\draw[thin] (5,-1)--(7,1)--(9,-1)--(11,1)--(13,-1);

\draw[fill=black] (5,-1) circle (1pt);
\draw[fill=black] (7,1) circle (1pt);
\draw[fill=black] (8,0) circle (1pt);
\draw[fill=black] (9,-1) circle (1pt);
\draw[fill=black] (11,1) circle (1pt);
\draw[fill=black] (12,0) circle (1pt);
\draw[fill=black] (13,-1) circle (1pt);

\node at (9,-2) {$\sigma_{3}(t) = N^{2}S^{2}N^{2}S^{2}$};
\end{tikzpicture}
\end{center}
\caption{The image under $\sigma_{3}$ of the 3-ary tree $t = \varepsilon \wedge \varepsilon( \varepsilon \wedge \varepsilon \wedge \varepsilon)$.}\label{fig:sigma}
\end{figure}
\end{example}

\begin{lemma}
The map $\sigma_{m} \colon \mathbf{B}^{m} \to \mathbf{D}^{m-1}$ sends $m$-ary trees with $n$ leaves to $(m-1)$-Dyck paths of length $n-1.$
\begin{proof}
We argue by induction. Recall that $n= m + g(m-1)$ for some integer $g \geq 0$. We  prove the result by induction on $g$. When $g=0$ there is only one tree to consider,  namely $t = \varepsilon \wedge \varepsilon \wedge \dots \wedge \varepsilon$.
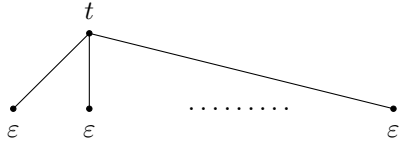
\begin{center}
\begin{tikzpicture}
\draw[fill=black] (0,0) circle (1pt);
\draw[fill=black] (-1,-1) circle (1pt);
\draw[fill=black] (0,-1) circle (1pt);
\draw[fill=black] (4,-1) circle (1pt);

\draw[thin] (0,0) -- (-1,-1);
\draw[thin] (0,0) -- (0,-1);
\draw[thin] (0,0) -- (4,-1);

\node at (0,0.3) {$t$};
\node at (-1,-1.3) {$\varepsilon$};
\node at (0,-1.3) {$\varepsilon$};
\node at (4,-1.3) {$\varepsilon$};
\node at (2,-1) {$\dots \dots \dots $};
\end{tikzpicture}
\end{center}
It is easy to see that $\sigma_{m}(t) = N^{m-1}S^{m-1}$ which is an $(m-1)$-Dyck path of length $m-1$. 

Now suppose that the result holds for $n=m+g^{\prime}(m-1)$ for all $g^{\prime} \leq g$. We consider the $g+1$ case. Let $t$ be an $m$-ary tree with $m+(g+1)(m-1)$ leaves. Then we may write $t = t_1 \wedge t_2 \wedge \dots \wedge t_m$ with the $t_i \in B_{n_i}^{m}$ and $n_1 + n_2 \dots + n_m = m+(g+1)(m-1).$ Then by definition $\sigma_{m}(t) = N^{m-1} \sigma_{m}(t_1) S \sigma_{m}(t_2) S \dots S \sigma_{m}(t_m)$ and by the inductive hypothesis, each $\sigma (t_i)$ is an $(m-1)$-Dyck paths of length $n_i -1$. In the expression for $\sigma_{m}(t)$ we have $m-1$ down-steps $S$ following the $N^{m-1}$ inbetween the $\sigma_{m} (t_i)$. Therefore the length of $\sigma_{m}(t)$ is $(n_1 -1)+(n_2 -1) + \dots (n_m -1) + (m-1) =m+(g+1)(m-1)-1$ as required. So $\sigma_{m}$ is indeed a map from $\mathbf{B}_{n}^{m}$ to $\mathbf{D}_{n-1}^{m-1}$.
\end{proof}
\end{lemma}

By the above lemma, $\sigma_{m}$ induces a map $\sigma_{m,n} \colon \mathbf{B}^{m}_{n} \to \mathbf{D}^{m-1}_{n-1}$. This map is in fact a bijection between $\mathbf{B}^{m}_{n}$ and $\mathbf{D}^{m}_{n-1}$. Let $t=t_1 \wedge t_2 \wedge \dots t_m$, where for $1 \leq i \leq m$ each $t_i$ is an $m$-ary tree with $n_i$ leaves and $n_1 + \dots + n_m = n$. Then $\sigma_{m,n}$ is defined as follows:

$$
\sigma_{m,n}(t) =
 N^{m-1} \sigma_{m,n_1}(t_1) S \sigma_{m,n_2}(t_2) S \dots S \sigma_{m,n_m}(t_m) 
.$$
\begin{proposition}\label{treePathBij}
The map $\sigma_{m,n} \colon \mathbf{B}^{m}_{n} \to \mathbf{D}^{m-1}_{n-1}$ is a bijection.
\begin{proof}
It is well known that both the finite sets $\mathbf{B}_{n}^{m}$ and $\mathbf{D}_{n-1}^{m-1}$ have cardinality $\frac{1}{(m-1)n+1}\binom{mn}{n}$, see for example \cite[Section 3]{heubach2008staircase}. So to show that $\sigma_{m,n}$ is a bijection it suffices to show that it is a surjection. We argue by induction on $n$. When $n=0$, it is trivial. 

Let $D \in \mathbf{D}_{n-1}^{m-1}$. Then the first step of $D$ is an up-step $N^{k_1(m-1)}$ where $k_1 \geq 1$ is an integer. So we can write $D= N^{m-1}N^{k_1(m-1) - (m-1)} \dots S$ as in Lemma \ref{form}. Let $(x^{\prime}_{1}, m-2)$ be the first point on $D$ with $y$-coordinate $m-2$ after the point $(m-1,m-1)$. Then the step in $D$ ending at $(x^{\prime}_{1}, m-2)$ must be a down-step $S$ starting at $(x_1,y_1)=(x^{\prime}_{1}+1,m-1)$. The part of $D$ from $(m-1,m-1)$ to $(x_1,y_1)$ is a translated $(m-1)$-Dyck path $D_1$. So we see that the path $D$ starts as $N^{m-1}D_{1}S$. Let $(x^{\prime}_{2}, m-3)$ be the first point on $D$ with $y$-coordinate $m-3$ after the point $(x^{\prime}_{1},m-2)$. As above, the step in $D$ ending at $(x^{\prime}_{2}, m-3)$ must be a down-step $S$ starting at $(x_2,y_2)=(x^{\prime}_{2}-1, m-2)$. The part of $D$ from $(x^{\prime}_{1},m-2)$ to $(x_2,y_2)$ is a translated $(m-1)$-Dyck path $D_2$. Therefore $D=N^{m-1}D_{1}SD_{2}S \dots S$, and continuing this argument we see that can be write $D=N^{m-1}D_{1}SD_{2}S \dots D_{m}S,$ where the $D_i$ are translated $(m-1)$-Dyck paths for $1 \leq i \leq m$.

Regarding the translated $(m-1)$-Dyck paths $D_i$ as $(m-1)$-Dyck paths, they each have length $n_{i} < n$ for $1 \leq i \leq m$. So by the inductive hypothesis, for each $D_i$ there exists an $m$-ary tree $t_i$ such that $D_{i} = \sigma_{m,n_{i}}(t_{i})$. It then follows that $D = \sigma_{m,n}(t_1 \wedge t_2 \wedge \dots \wedge t_m)$.

\end{proof}
\end{proposition}
Going forward, we drop the subscripts on $\sigma_{m,n}$ and just write $\sigma$ when it clear what is meant from the context. 
\begin{proposition}{\label{depthExp}}
Let $t$ be an $m$-ary tree with $n$ leaves and depth $( \delta^{l_1}(t), \delta^{l_2}(t) \dots \delta^{l_m}(t) )$. Then $\sigma (t) = N^{d_1}SN^{d_2}S \dots SN^{d_{n-1}}SN^{d_n}$ where the $d_i$ are given by  

$$ d_1 = (m-1)\delta^{l_1}_{1}(t),$$
$$d_j =  \left(\sum\limits_{i=1}^{m} (m-i)(\delta^{l_{i}}_{j}(t) - \delta^{l_i}_{j-1}(t))\right)+1, \text{ for } 2 \leq j \leq n.$$
\begin{proof}
Recall that $n$ satisfies the equation $n = m + g(m-1)$ for some integer $g \geq 0$. We prove the result by induction on $g$. When $g=0$ there is only one tree to consider,  namely $t = \varepsilon \wedge \varepsilon \wedge \dots \wedge \varepsilon$.
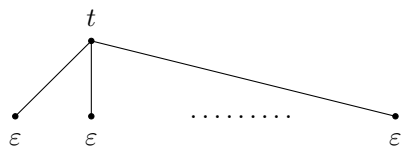
\begin{center}
\begin{tikzpicture}
\draw[fill=black] (0,0) circle (1pt);
\draw[fill=black] (-1,-1) circle (1pt);
\draw[fill=black] (0,-1) circle (1pt);
\draw[fill=black] (4,-1) circle (1pt);

\draw[thin] (0,0) -- (-1,-1);
\draw[thin] (0,0) -- (0,-1);
\draw[thin] (0,0) -- (4,-1);

\node at (0,0.3) {$t$};
\node at (-1,-1.3) {$\varepsilon$};
\node at (0,-1.3) {$\varepsilon$};
\node at (4,-1.3) {$\varepsilon$};
\node at (2,-1) {$\dots \dots \dots $};
\end{tikzpicture}
\end{center}
For this tree we have that $\delta^{l_i}_{j} =\delta_{ij}$, where the right hand side is the usual Kronecker delta function. We also have that $\sigma(t) = N^{m-1}SN^{0}SN^{0}S \dots N^{0}SN^{0} = N^{m-1}S^{m-1}$. We now need to verify that the exponents of the $N$s satisfy the above relations. Indeed $d_1= m-1 = (m-1)\delta^{l_1}_{1}$. Moreover $\sum\limits_{i=1}^{m} (m-i)(\delta^{l_{i}}_{j} - \delta^{l_i}_{j-1})+1 = (m-j) - (m-(j-1))+1 = (m-m) + ((j-1)-j)+1 = 0 = d_j$ for $2 \leq j \leq n$.

Now suppose that the result holds for $n=m+g^{\prime}(m-1)$ for all $g^{\prime} \leq g$. We consider the $g+1$ case. Let $t$ be an $m$-ary tree with $n= m+(g+1)(m-1)$ leaves. Then we may write $t =  t_1 \wedge \dots \wedge t_m$ where each $t_i$ is the subtree rooted at the $i^{\text{th}}$ child of the root of $t$. Each subtree $t_i$ has $n_i < n$ leaves and $n_1 + n_2 \dots + n_m = n$. In writing $t$ as the meet of its sub-trees at the root, we partition the leaves of $t$. We identify each leaf of $t$ with a pair $(h,j)$ if it lies in the subtree $t_h$ and it is the $j^{\text{th}}$ leaf in the linear order on the leaves of $t_h$ where $1 \leq j \leq n_{h}$. Therefore for the leaf $(h,j)$ we have that, 
\begin{equation}\label{eqn:depthcases}
\delta^{l_i}_{h,j}(t) =
\begin{cases}
\delta^{l_i}_{h,j}(t_h) +1  & \text{ if } i=h  ;\\
\delta^{l_i}_{h,j}(t_h) & \text{otherwise}.\\
\end{cases}
\end{equation}
By the inductive hypothesis $\sigma(t_h) = N^{d_{h,1}}SN^{d_{h,2}} \dots SN^{h,n_h}$, where $$d_{h,1}(t) =(m-1) \delta^{l_1}_{h,1}(t_h)$$ and 
$$d_{h,j}(t) = \sum\limits_{i=1}^{m} (m-i)(\delta^{l_{i}}_{h,j}(t_h) - \delta^{l_i}_{h,j-1}(t_h))+1, \text { for } 2 \leq j \leq n_{h}.$$
By definition, $\sigma(t)=N^{m-1}\sigma(t_1)S\sigma(t_2)\dots S\sigma(t_m)$, so
 $$\sigma(t) = N^{m-1}N^{d_{1,1}}SN^{d_{1,2}}S \dots SN^{d_{1,n_1}}SN^{d_{2,1}}SN^{d_{2,2}}S \dots N^{d_{2,n_{2}}}SN^{d_{3,1}}S \dots SN^{d_{m,n_m}}$$
$$ \hspace{0.75cm} = N^{(m-1)+d_{1,1}}SN^{d_{1,2}}S \dots SN^{d_{1,n_1}}SN^{d_{2,1}}SN^{d_{2,2}}S \dots N^{d_{2,n_{2}}}SN^{d_{3,1}}S \dots SN^{d_{m,n_m}}$$
 Now we verify that the exponents of the $N$s satisfy the required relations. We see that
$$(m-1) + d_{1,1} = (m-1) +(m-1) \delta^{l_1}_{1,1}(t_1) = (m-1)(\delta^{l_1}_{1,1}(t_1)+1) = (m-1)\delta^{l_1}_{1,1}(t),$$
so the first exponent satisfies the required relation. We also see that
$$d_{h,j} = \sum\limits_{i=1}^{m} (m-i)(\delta^{l_{i}}_{h,j}(t_h) - \delta^{l_i}_{h,{(j-1)}}(t_h))+1 = \sum\limits_{i=1}^{m} (m-i)(\delta^{l_{i}}_{h,j}(t) - \delta^{l_i}_{h,{(j-1)}}(t))+1, \text { for } 2 \leq j \leq n_{h}, \\$$
by (\ref{eqn:depthcases}), therefore the $d_{j,h}$ also satisfy the required relation for $t$. 

The only exponents left to verify are the $d_{h,1}$ for $2 \leq h \leq m$. In this case, the leaf $(h-1,n_{h-1})$ is the rightmost leaf in the subtree $t_{h-1}$, so by Lemma \ref{rem1}, $\delta^{l_i}_{h-1,n_{h-1}}(t_{h-1}) = 0$ when $i\neq m$. Therefore by (\ref{eqn:depthcases}), $\delta^{l_i}_{h-1,n_{h-1}}(t) = 0$ when $i\neq m,h-1$, so $\delta^{l_{h-1}}_{h-1,n_{h-1}}(t) = 1$. The leaf $(h,1)$ is the leftmost leaf in the subtree $t_h$, so by a dual statement of Lemma \ref{rem1}, $\delta^{l_i}_{h,1}(t_h) = 0$ when $i \neq 1$. Therefore by (\ref{eqn:depthcases}),  $\delta^{l_i}_{h,1}(t) = 0$ when $i \neq 1,h$ and $\delta^{l_{h}}_{h,1}(t)=1$. So it follows that, 

\begin{align*}
\sum\limits_{i=1}^{m} (m-i)(\delta^{l_i}_{h,1}(t)-\delta^{l_i}_{h-1,n_{h-1}}(t)) + 1 &= (m-1)\delta^{l_1}_{h,1}(t_h)+(m-h)-(m-(h-1)) \\
&\hspace{0.5cm} - (m-m)\delta^{l_m}_{h-1,n_{h-1}}(t_{h-1})+1 \\
 &= (m-1)\delta^{l_1}_{h,1}(t_h) \\
&= d_{h,1}
\end{align*}
therefore the $d_{h,1}$ also satisfy the required relations for $2 \leq h \leq m$. 
This completes the proof. 
\end{proof}
\begin{remark}\label{dn=0}
It is important to note that $d_n = 0$ in Proposition \ref{depthExp} since otherwise $D$ is not a Dyck path. We can observe that $d_n=0$ by referencing Lemma \ref{rem1} and Lemma \ref{rem2}. So in the above proposition, $\sigma (t)$ is indeed a Dyck path. Also note that since $d_n=0$, this form of $\sigma(t)$ is the same as that given in Lemma \ref{form}. 
\end{remark}
\end{proposition}
The bijection between $m$-ary trees with $n$ leaves and $(m-1)$-Dyck paths of length $n-1$ induces an operation corresponding to $k$-rotation on Dyck paths, which we shall call a \textit{$k$-compression}. Recall in the definition of a right $k$-rotation we replace a sub-tree of the form, $$s = t_1 \wedge t_2 \wedge t_3 \wedge\dots\wedge t_{j-1} \wedge(t_j \wedge t_{j+1} \wedge\dots\wedge t_{j+k(m-1)}) \wedge t_{j+k(m-1)+1} \wedge \dots \wedge t_{(m-1)+k(m-1)}$$ by a subtree of the form,
$$ s^{\prime} = t_1 \wedge t_2 \wedge\dots\wedge t_{j} \wedge (t_{j+1} \wedge t_{j+2}\dots\wedge t_{j+k(m-1)+1}) \wedge t_{j+k(m-1)+2} \wedge\dots\wedge t_{(m-1)+k(m-1)}
.$$
It is easy to see that $$\sigma(s) = N^{m-1}D_1SD_2S \dots D_{j-1}SN^{k(m-1)}D_{j}SD_{j+1}S \dots SD_{m+k(m-1)},$$ and 
$$\sigma(s^{\prime}) = N^{m-1}D_1SD_2S \dots D_{j-1}SD_{j}SN^{k(m-1)}D_{j+1}S \dots SD_{m+k(m-1)},$$ where $D_i =\sigma(t_i)$. 

\begin{definition}{\textbf{Right $k$-Compression}}. 
Let $k \geq 1$ and $1 \leq j \leq m-1$ be positive integers. Let $D$ be an $(m-1)$-Dyck path of length $n-1$. Suppose $D$ contains a sub-Dyck path of the form $$X = N^{m-1}D_1SD_2S \dots D_{j-1}SN^{k(m-1)}D_{j}SD_{j+1}S \dots SD_{m+k(m-1)},$$ where the $D_{i}$ are (possibly translated) Dyck paths which may be empty. Then a \textit{right $k$-compression} at $X$ is the operation of replacing $X$ with the sub-Dyck path
 $$X^{\prime}  = N^{m-1}D_1SD_2S \dots D_{j-1}SD_{j}SN^{k(m-1)}D_{j+1}S \dots SD_{m+k(m-1)}.$$ 
\end{definition}
The inverse operation of replacing $X^{\prime}$ with $X$ will be called a \textit{left $k$-compression}. Let $D, D^{\prime}$ be $(m-1)$-Dyck paths of length $n-1$. Write $D \preceq_{k} D^{\prime}$ to mean that $D^{\prime}$ can be obtained from $D$ by applying finitely many right $k$-compressions. The induced partial order on $\mathbf{D}^{m-1}_{n-1}$ is called the $k$-associative order. The connected components of $\mathbf{D}^{m-1}_{n-1}$ under the $k$-associative order are called $k$-components. We will say two $(m-1)$-Dyck paths are $k$-equivalent if they belong to the same $k$-component. 

Let $\mathbb{M} \subset \mathbb{N}^{n}$ be the set of $n$-tuples of non-negative integers $(e_1, e_2, \dots, e_n)$ satisfying the following relations, 
$$e_1 + e_2 + \dots + e_n =n-1,$$
$$(m-1) | e_i \text{ for } 1 \leq i \leq n,$$
$$e_1 + e_2 + \dots + e_{j-1} \geq j-1 \text{ for all } 1 \leq j \leq n.$$ 
Notice that it follows from the first and last relation that $e_n=0$.
\begin{proposition}{\label{degreemap}}
The map $d \colon \mathbf{D}^{m-1}_{n-1} \to \mathbb{M}$ maps an $(m-1)$-Dyck path of length $(n-1)$ $D = N^{d_1}SN^{d_2}S \dots SN^{d_n}$ to the $n$-tuple $d(D)= (d_1, \dots, d_n)$. This map is a bijection.
\begin{proof}
Let $D= N^{d_1}SN^{d_2}S \dots SN^{d_n}$ be an $(m-1)$-Dyck path. Let $d(D)=(d_1,d_2, \dots, d_n)$. Note that $d_{n}=0$ by Remark \ref{dn=0}, so the form of $D$ is precisely as in Lemma \ref{form}. By Lemma \ref{form} the tuple $(d_1, d_2, \dots, d_{n-1})$ is unique, so the map $d$ is well-defined. Since $D$ is an $(m-1)$-Dyck path, by definition $(m-1) | d_i$. All Dyck paths start and end on the $x$-axis, therefore they must go up the same number of times as they go down. Hence if a path has length $n-1$, which is the number of down-steps $S$, then $ d_1 + \dots + d_n = n-1$. By definition, Dyck paths cannot go below the $x$-axis, this is to say that $d_1 + \dots + d_{j-1} \geq j-1$ for all $j \geq 1$.

Let $f \colon \mathbb{M} \to \mathbf{D}^{m-1}_{n-1}$ be the map given by $f(e_1,e_2,\dots,e_n)= N^{e_1}SN^{e_2}S \dots SN^{e_n}.$ This is a $(m-1)$-dyck path by the arguments similar to those above. It is easy to see that $f(d(D)) = D$, and $d(f((e_1, \dots , e_n))) = (e_1, \dots , e_n)$. Therefore $d$ is indeed a bijection.
\end{proof}
\end{proposition}

\begin{proposition}{\label{imp}}
Let $D , D^{\prime}$ be $(m-1)$-Dyck paths of length $n-1$ with $d(D) = (d_1, \dots , d_n)$ and $d(D^{\prime}) = (d_1^{\prime}, \dots , d_n^{\prime})$. Suppose that we can obtain $D^{\prime}$ from $D$ by applying a right $k$-compression to $D$. Then there exist $1\leq j < i \leq n$ such that $d_{i}^{\prime} = d_i +k(m-1)$, $d_{j}^{\prime} = d_j - k(m-1)$ and $d_{h}^{\prime} = d_h$ for $h \neq i,j.$  
\begin{proof}
Recall that in the definition of right $k$-compression, we replace $$X= N^{m-1}D_1SD_2S \dots D_{a-1}SN^{k(m-1)}D_{a}SD_{a+1}S \dots SD_{m+k(m-1)}$$ with $$X^{\prime} = N^{m-1}D_1SD_2S \dots D_{a-1}SD_{a}SN^{k(m-1)}D_{a+1}S \dots SD_{m+k(m-1)},$$ thus moving the substring $N^{k(m-1)}$ from the immediate left of the (possibly translated) Dyck path $D_a$ to the immediate left of (possibly translated) Dyck path $D_{a+1}$. Let $D = N^{d_1}SN^{d_2}S \dots SN^{d_n}.$ We have the sub-strings $N^{k(m-1)}D_a = N^{d_j}S \dots S$ and $D_{a+1} =N^{d_i} \dots S$ in $D$. So in moving $N^{k(m-1)}$ we get the sub-strings $D_a = N^{d_j-k(m-1)}S \dots S$ and $N^{k(m-1)}D_{a+1} =N^{d_i + k(m-1)}S \dots S$ in the Dyck path $D^{\prime}$. This proves the statement of the proposition. 
\end{proof}
\end{proposition}
Remark that if we replace right $k$-compression with left $k$-compression in the above proposition, we have that $j > i$ instead.  

\begin{corollary}
Let $D , D^{\prime}$ be $(m-1)$-Dyck paths of length $n-1$. If $D$ and $D^{\prime}$ are $k$-equivalent then $d(D) \equiv d(D^{\prime})$ mod $k(m-1)$. 
\begin{proof}
This is an immediate consequence of Proposition \ref{imp}.
\end{proof}
\end{corollary}

Let $D$ be an $m$-Dyck path of length $n$. We say that $D$ is \textit{$k$-minimal} if it is minimal in its $k$-equivalence class. That is to say there does not exist a Dyck path $D^{\prime} \in \mathbf{D}_{n}^{m}$ such that $D^{\prime} \preceq_{k} D$. Let $p=(x,y)$ in  $\mathbb{Z}^{2}$ be a point on the $m$-Dyck path $D$. The \textit{level} of the point $p$ is the integer $y$, and we say that $p$ is on the $y^{\text{th}}$ level. 

\begin{proposition}{\label{minimal}}
An $(m-1)$-Dyck path D is minimal if and only if for $d(D) = (d_1, \dots , d_n),$ we have that $d_i < k(m-1) \text{ for all } i \neq 1.$
\begin{proof}
Suppose that we have that $d_i < k(m-1)$ for all $i \neq 1$ and $D$ is not minimal. Then we can left $k$-compress $D$ to obtain another dyck path $D^{\prime}$. By Proposition \ref{imp} there is some $j > 1$ such that the $j$-th entry of $d(D^{\prime})$ is $d_{j}^{\prime} = d_{j} -k(m-1)$. By the assumption that  $d_i < k(m-1)$ for $i \neq 1$, we must have that $d_{j}^{\prime} < 0$, a contradiction. Therefore $D$ must be minimal.

Recall $D$ is of the form $D=N^{d_1}S \dots SN^{d_{i}}S\dots SN^{d_n}$.  Suppose that $D$ is minimal and there exists some $i \neq 0$ such that $d_i \geq k(m-1)$. We will show that $D$ is not minimal by demonstrating that we can left $k$-compress $D$. That is to say we will show that we have a sub-Dyck path $X^{\prime}$ required to perform a left $k$-compression, where 
 $$X^{\prime} = N^{m-1}D_{1}SD_{2}S \dots D_{j-1}SD_{j}SN^{k(m-1)}D_{j+1}S \dots SD_{(m-1) +k(m-1)}$$ for $1 \leq j \leq (m-1)$.
 

Suppose the up-step $N^{d_i}$ starts at some point $(b,l)$ and ends at $(b+d_i, l +d_i)$. The immediately preceding down-step $S$ starts at $(b-1,l+1)$ and ends at $(b,l)$. Let $0 \leq x \leq b-1$ be maximal such that the point $(x,l)$ is on the Dyck path $D$. By the maximality, the point $(x,l)$ is part of an up-step.  Let $U$ to be the up-step in $D$ beginning at $(x,l)$ if $(x,l)$ is at the start of an up-step; otherwise let $U$ to be the up-step containing $(x,l)$. Let $(x_1,y_1)$ be the end point of the up-step $U$. Let $(x_0,y_0) = (x_1 -(m-1),y_1 -(m-1))$, this is the start point of the up-step $U$. See the figure below.
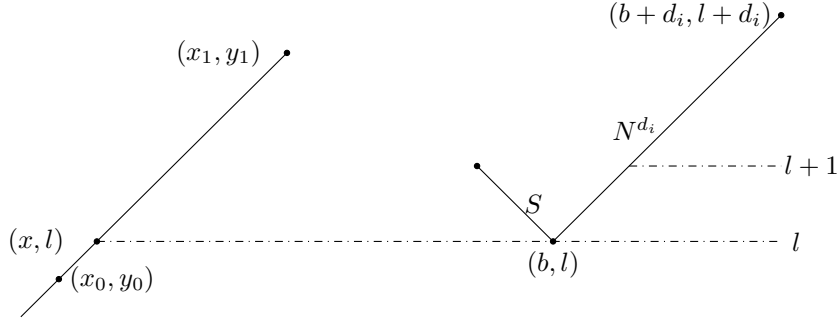
\begin{center}
\begin{tikzpicture}
  
\draw[fill=black] (0,0) circle (1pt);
\draw[fill=black] (2.5,2.5) circle (1pt);  
\draw[fill=black] (5,1) circle (1pt);
\draw[fill=black] (6,0) circle (1pt);
\draw[fill=black] (9,3) circle (1pt);
\draw[fill=black] (-0.5,-0.5) circle (1pt);

\node at (1.6,2.5) {$(x_1,y_1)$};  

  \node at (-0.8,0) {$(x,l)$};
  \node at (9.2,0) {$l$};
  \node at (9.4,1) {$l+1$};
  \node at (6,-0.3) {$(b,l)$};
  \node at (7.8,3) {$(b+d_i, l +d_i)$};
  \node at (0.2,-0.5) {$(x_0,y_0)$};
  
 \draw[thin] (-1,-1)--(0,0);
  \draw[thin] (0,0)--(2.5,2.5);
  \draw[thin] (5,1) --(6,0)node[midway, right] {$S$};
  \draw[thin](6,0) --(9,3)node[midway, left] {$N^{d_i}$};
  \draw[thin, dash dot] (0,0)--(9,0);
  \draw[thin, dash dot] (7,1)--(9,1);
\end{tikzpicture} 
\end{center}


Let $(x_2,y_2)$ be the point at which it is the first time the Dyck path goes below the level $y_1$ after the point $(x_1,y_1)$. That is $y_2 = y_{1}-1$. Then the subpath $D_1$ starting from $(x_1,y_1)$ and ending at $(x_2-1,y_2+1)$ is a translated $(m-1)$-Dyck path. Note that it could happen that $(x_1,y_1) = (x_2-1,y_2+1)$, in this case $D_1$ is just the empty $(m-1)$-Dyck path.


Let $(x_3,y_3)$ be the point at which it is the first time the Dyck path goes below the level $y_{1}-1$ after the point $(x_2,y_2)$, that is $y_3 = y_{1}-2$. We define $D_2$ to be the path starting from $(x_2,y_2)$ to $(x_3-1, y_3+1)$. As before $D_2$ is a translated $m$-Dyck path which starts and ends on the $(y_{1}-1)^{\text{th}}$ level. 
 
Let $j = y_{1}-l$. We can repeat this procedure to define translated $m$-Dyck paths $D_{3}, D_{4}, \dots, D_{j}$. Here each $m$-Dyck path $D_r$ starts at the point $(x_r,y_r)$ and ends at the point $(x_{r+1}-1,y_{r+1}+1)$, where the start and end points are defined as above and $y_{r+1} = y_{r} -1 = y_{1} -r$. Note that the translated $m$-Dyck path $D_j$ begins on the level $y_{j} = y_{1}-(j-1)=l+1$, so the last point of $D_j$ is $(x_{j+1}-1,l+1)$ for some $x_{j+1}-1 \leq b$. We claim that $(x_{j+1}-1,l+1)=(b-1,l+1)$. The point $(b-1,l+1)$ is the last time we are on the $(l+1)^{\text{th}}$ level before the $N^{d_i}$ up-step. By construction, there is a down-step from $(x_{j+1}-1,l+1)$ to $(x_{j+1},l)$. By the maximality of $x$ we have that $x_{j+1} = b$ or $x_{j+1} = x $. Note that $x_{j+1} \geq x_{1} > x$, so we have that $x_{j+1} = b$.

 
So far we have constructed a subpath  from $(x_0,y_0)$ to $(b,l+1)$ given by 
$$X^{\prime \prime} = N^{m-1}D_{1}SD_{2}S \dots SD_{j}S,$$ 
where the $S$ down-steps are the down steps from $(x_i-1,y_i+1)$ to  $(x_i,y_i)$. Note that $y_i = y_{1} -(i-1)$ for $2 \leq i < j$ and the $S$ after $D_j$ is the one from $(b-1,l+1) \text{ to } (b,l)$. The $N^{m-1}$ is the up-step $U$ from $(x_0,y_0)$ to $(x_1,y_1)$.

Since $d_i \geq k(m-1)$, there is an up-step $N^{k(m-1)}$ from $(b,l)$ to $(x_{j+1},y_{j+1}) =(b+k(m-1), l+k(m-1))$. We define $D_{j+1}$ to be the path from $(x_{j+1},y_{j+1})$ to $(x_{j+2}-1,y_{j+2}+1)$ where $(x_{j+2},y_{j+2})$ is the point at which the Dyck path first sits on level $l+k(m-1)-1$ after $(x_{j+1}, y_{j+1})$. In the same fashion we define the $(m-1)-j+k(m-1)$ sub paths $D_{j+2}, D_{j+3} \dots D_{(m-1)+k(m-1)}$. These are all translated $m$-Dyck paths by the same arguments as above. By how we construct the Dyck paths, we see that the path $D_{(m-1)+k(m-1)}$ ends on level $y_0 = y_1 - (m-1) $.

So we have successfully constructed the sub-Dyck path of $D$, $$X^{\prime} = N^{m-1}D_{1}SD_{2}S \dots D_{j-1}SD_{j}SN^{k(m-1)}D_{j+1}S \dots SD_{(m-1) +k(m-1)}.$$ As before the $S$ are the intermediate down steps between the $D_i$ and the $D_i$ may also be empty.
 \end{proof}
We illustrate the constructive proof above with an example for the case where $k=2$ and $m=3$. 
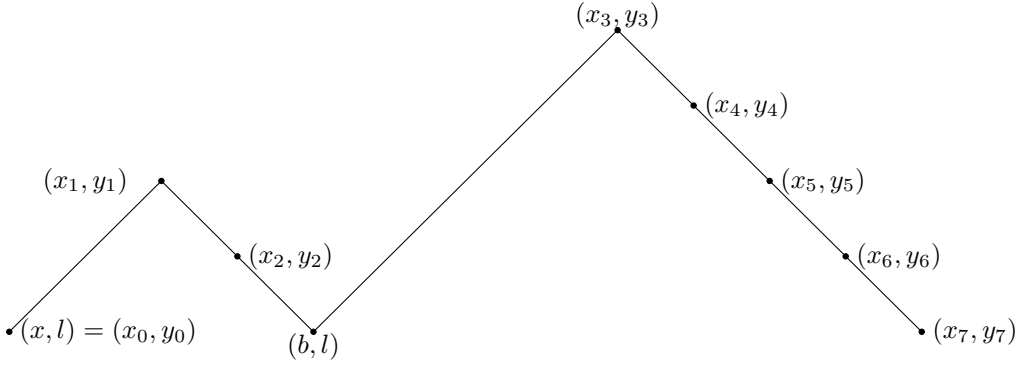
\begin{center}
 \begin{tikzpicture}
\draw[fill=black] (0,0) circle (1pt);
\draw[fill=black] (2,2) circle (1pt);  
\draw[fill=black] (3,1) circle (1pt);  
\draw[fill=black] (4,0) circle (1pt);  
\draw[fill=black] (8,4) circle (1pt);  
\draw[fill=black] (9,3) circle (1pt);  
\draw[fill=black] (10,2) circle (1pt);  
\draw[fill=black] (11,1) circle (1pt);  
\draw[fill=black] (12,0) circle (1pt);  

  \node  at (1,2) {$(x_1,y_1)$};
   \node at (1.3,0) {$(x,l)=(x_0,y_0)$};
   \node at (3.7,1) {$(x_2,y_2)$};
   \node at (8,4.2) {$(x_3,y_3)$}; 
   \node at (9.7,3) {$(x_4,y_4)$};
    \node at (10.7,2) {$(x_5,y_5)$};
     \node at (11.7,1) {$(x_6,y_6)$};
      \node at (12.7,0) {$(x_7,y_7)$};
      \node at (4,-0.2) {$(b,l)$};

 \draw[thin] (0,-0)--(2,2)--(4,0)--(8,4)--(12,0);
 
\end{tikzpicture}
\end{center}
The $N^{m-1} = N^{2}$ up-step is the one from $(x_0,y_0)$ to $(x_1,y_1)$. The translated $2$-Dyck paths $D_1$ and $D_2$ are the empty paths $N^{0}S^{0}$ at $(x_1,y_1)$ and $(x_2,y_2)$ respectively. The $N^{k(m-1)}=N^{2(2)}$ up-step is the one from $(b,l)$ to $(x_3,y_3)$. The rest of the translated $2$-Dyck paths $D_3,\dots,D_6$ are the empty paths at $(x_4,y_4), \dots ,(x_7,y_7)$ respectively. Therefore, $X^{\prime}$ in this case is the whole path above.  
\end{proposition}

We now show that minimal $m$-Dyck paths do exist and that they are unique in each $k$-equivalence class. 

\begin{proposition}{\label{minimalunique}}
Every $k$-equivalence class contains a unique minimal Dyck path.
\begin{proof}
To show existence, we consider the bijection $d \colon \mathbf{D}_{n-1}^{m-1} \to \mathbb{M}$ from Proposition \ref{degreemap}  where $d(D) = (d_1, \dots ,d_n)$. Endow $\mathbb{N}^{n}$ with the standard lexicographic order. Then by proposition \ref{imp}, $d$ is order reversing. That is $D \prec_{k} D^{\prime}$ implies that $d(D^{\prime}) <_{\text{lex}} d(D)$. Recall that the lexicographic order is a partial order therefore it has no cycles because of anti-symmetry. So suppose that there is a $k$-equivalence class with no minimal Dyck path. Take $D$ belonging to such a class and repeatedly left $k$-compress it. Since there is no minimal element in this class, we can do this indefinitely. As a result we obtain the descending chain, $ \dots \prec_{k}  D^{(a)} \prec_{k} \dots \prec_{k} D^{1} \prec_{k} D$. Applying $d$ to the descending chain, we get the ascending chain $d(D) <_{\text{lex}} d(D^{1}) <_{\text{lex}} \dots <_{\text{lex}} d(D^{a}) <_{\text{lex}} \dots$. Since $\mathbf{D}_{n-1}^{m-1}$ is finite, this ascending chain must then be a cycle. Since the lexicographic order is anti-symmetric, this cycle must contain only one element. Therefore if a $k$-equivalence class doesn't have a minimal element, it only contains one Dyck path in which case that Dyck path is trivially minimal. This is a contradiction to our assumption. Thus every $k$-equivalence class has a minimal Dyck path. 

Suppose we have two minimal Dyck paths $D$ and $D^{\prime}$ in an equivalence class. By Proposition \ref{minimal} all but the first entries of $d(D)$ and $d(D^{\prime})$ are strictly less than $k(m-1)$. But since $D$ and $D^{\prime}$ are $k$-equivalent, $d(D) \equiv d(D^{\prime})$ mod $k(m-1)$. This means all but the first entries of $d(D)$ and $d(D^{\prime})$ are equal. The equality of these entries forces the first entries to also be equal since clearly it cannot be the case otherwise. Therefore $d(D) = d(D^{\prime})$ which implies $D=D^{\prime}$. So the minimal Dyck paths are unique in their equivalence classes. 
\end{proof}
\end{proposition}

\begin{theorem}{\label{DyckChar}}
Let $D, D^{\prime} \in \mathbf{D}_{n}^{m}$  be two Dyck paths. Then $D$ and $D^{\prime}$ are $k$-equivalent if and only if $d(D) \equiv d(D^{\prime}) \text{ mod } k(m-1)$. 

\begin{proof}
Suppose that $D$ and $D^{\prime}$ are $k$-equivalent. Then suppose without loss of generality that we obtain $D^{\prime}$ from $D$ by application of a finite sequence of $k$-compressions. From Proposition \ref{imp}, we see that a $k$-compression maps $d(D)$ to an $n$-tuple which is congruent to $d(D)$ modulo $k(m-1)$. Therefore since $d(D^{\prime})$ an  $n$-tuple which is a result of a finite sequence of $k$-compressesions on $D$, then $d(D) \equiv d(D^{\prime}) \text{ mod } k(m-1)$.

Suppose now that $d(D) \equiv d(D^{\prime})$ mod $k(m-1)$. Consider their respective minimal representatives in their $k$-equivalence classes $D_{min}$ and $D_{min}^{\prime}$ respectively. Then $d(D_{min}) \equiv d(D) \equiv d(D^{\prime}) \equiv d(D_{min}^{\prime})$ mod $k(m-1)$. Therefore $d(D_{min}) \equiv d(D_{min}^{\prime})$ mod $k(m-1)$, hence by Proposition \ref{minimal} we obtain that $d(D_{min}) = d(D_{min}^{\prime})$ which means $D_{min} = D_{min}^{\prime}$. Therefore $D$ and $D^{\prime}$ belong to the same $k$-equivalence class. So we conclude $D$ and  $D^{\prime}$ are $k$-equivalent. 
\end{proof}
\end{theorem}

\begin{theorem}{\label{TreeChar}}
Let $t,t^{\prime}$ be a pair of $m$-ary trees with $n$ leaves and depth $( \delta^{l_1}(t), \delta^{l_2}(t),\dots,\delta^{l_{m}}(t) )$, $(\delta^{l_1}(t^{\prime}), \delta^{l_2}(t^{\prime}),\dots,\delta^{l_{m}}(t^{\prime}))$ respectively. Then $t$ and $ t^{\prime}$ are $k$-equivalent if and only if \[  \sum\limits_{i=1}^{m-1} (m-i)\delta^{l_{i}}(t) \equiv \sum\limits_{i=1}^{m-1} (m-i)\delta^{l_{i}}(t^{\prime}) \text{ mod } k(m-1),  \] where the addition on the $n$-tuples is componentwise. 
\begin{proof}
Suppose that $t$ and $t^{\prime}$ are $k$-equivalent, then their corresponding Dyck paths $D=\sigma(t)$ and $D^{\prime} = \sigma(t^{\prime})$ respectively are also $k$-equivalent. Therefore by Theorem \ref{DyckChar} $d(D) \equiv d(D^{\prime})$ mod $k(m-1)$. 
By Proposition \ref{depthExp},
$$ d_1 = (m-1)\delta^{l_1}_{1}(t),$$
$$d_j =  \sum\limits_{i=1}^{m} (m-i)(\delta^{l_{i}}_{j}(t) - \delta^{l_i}_{j-1}(t))+1,\text{ for } j >1.$$
Since $d_{1} \equiv d_{1}^{\prime}$ mod $k(m-1)$,  we have that $(m-1)\delta^{l_1}_{1}(t) \equiv (m-1)\delta^{l_1}_{1}(t^{\prime})$ mod $k(m-1)$. 
Furthermore, we observe that from the structure of the of $m$-ary trees that $\delta^{l_i}_{1}(t) = 0$  and $\delta^{l_i}_{1}(t^{\prime}) = 0$ for $i \neq 1$. So we have that
$$
(m-1)\delta^{l_i}_{1}(t) \equiv (m-1)\delta^{l_i}_{1}(t^{\prime}) \text { mod }k(m-1) \text{ for } 1 \leq i \leq m,
$$
hence 
$$
\sum\limits_{i=1}^{m} (m-1)\delta^{l_i}_{1}(t) \equiv \sum\limits_{i=1}^{m} (m-1)\delta^{l_i}_{1}(t^{\prime}) \text { mod }k(m-1).
$$
From the fact that,\\
$d_{2} = \sum\limits_{i=1}^{m} (m-i)(\delta^{l_{i}}_{2}(t) - \delta^{l_i}_{1}(t))+1 \equiv d_{2}^{\prime} = \sum\limits_{i=1}^{m} (m-i)(\delta^{l_{i}}_{2}(t^{\prime}) - \delta^{l_i}_{1}(t^{\prime}))+1$ mod $k(m-1)$, 

we conclude that, 
$$\sum\limits_{i=1}^{m} (m-i)\delta^{l_{i}}_{2}(t) \equiv \sum\limits_{i=1}^{m} (m-i)\delta^{l_{i}}_{2}(t^{\prime}) \text{ mod } k(m-1).$$
From this congruence and the congruence $$d_3 =  \sum\limits_{i=1}^{m} (m-i)(\delta^{l_{i}}_{3}(t) - \delta^{l_i}_{2}(t))+1 \equiv d_{3}^{\prime} = d_3 =  \sum\limits_{i=1}^{m} (m-i)(\delta^{l_{i}}_{3}(t^{\prime}) - \delta^{l_i}_{2}(t^{\prime}))+1,$$
we conclude that, 
$$\sum\limits_{i=1}^{m} (m-i)\delta^{l_{i}}_{3}(t) \equiv \sum\limits_{i=1}^{m} (m-i)\delta^{l_{i}}_{3}(t^{\prime}) \text{ mod } k(m-1).$$
Continuing in this manner we obtain the following, 
$$\sum\limits_{i=1}^{m} (m-i)\delta^{l_{i}}_{j}(t) \equiv \sum\limits_{i=1}^{m} (m-i)\delta^{l_{i}}_{j}(t^{\prime}) \text{ mod } k(m-1),\text{ for }1 \leq j \leq n.$$

This is the same as saying, 
\[  \sum\limits_{i=1}^{m-1} (m-i)\delta^{l_{i}}(t) \equiv \sum\limits_{i=1}^{m-1} (m-i)\delta^{l_{i}}(t^{\prime}) \text{ mod } k(m-1).  \]

Now for the converse, suppose that \[  \sum\limits_{i=1}^{m-1} (m-i)\delta^{l_{i}}(t) \equiv \sum\limits_{i=1}^{m-1} (m-i)\delta^{l_{i}}(t^{\prime}) \text{ mod } k(m-1).  \]
This implies that, 
$$\sum\limits_{i=1}^{m} (m-i)\delta^{l_{i}}_{j}(t) \equiv \sum\limits_{i=1}^{m} (m-i)\delta^{l_{i}}_{j}(t^{\prime}) \text{ mod } k(m-1), \text{ for }1 \leq j \leq n.$$
This further implies that
$$ d_1 = (m-1)\delta^{l_1}_{1}(t) \equiv d_{1}^{\prime}  =  (m-1)\delta^{l_1}_{1}(t) \text{ mod }  k(m-1)$$
$$d_j =  \sum\limits_{i=1}^{m} (m-i)(\delta^{l_{i}}_{j}(t) - \delta^{l_i}_{j-1}(t))+1 \equiv d_j^{\prime} =  \sum\limits_{i=1}^{m} (m-i)(\delta^{l_{i}}_{j}(t^{\prime}) - \delta^{l_i}_{j-1}(t^{\prime}))+1 \text{ mod } k(m-1).$$
Therefore by Theorem \ref{DyckChar} we have that $D=\sigma(t) $ and $D^{\prime}=\sigma(t^{\prime})$ are $k$-equivalent which implies that $t$ and $t^{\prime}$ are $k$-equivalent. 
\end{proof}
\end{theorem}

\section{An Application to $m$-ary operations}
In this section we will introduce a particular $k$-associative $m$-ary operation which will be  denoted by $\circ$. This operation will be used to evaluate $m$-ary parenthesizations and we will show that this operation characterises $k$-equivalence. This is to say that two parenthesizations will be $k$-equivalent ($k$-associative) if and only if their evaluations under this operation are equal. For the rest of this section, we fix integers $m \geq 2$, $g \geq 0$, $k \geq 1$ and $n = m+g(m-1)$.

Let $A = \mathbb{C} \langle u_1,u_2,\dots,u_n \rangle$ be the free unital associative algebra over $\mathbb{C}$ in $n$ indeterminates $u_1,u_2,\dots,u_n$. We define a binary operation $\circ$ on $A$ as follows. Let $\omega$ be an element of $A$ of order $k(m-1)$ e.g. $\omega = e^{\frac{2\pi i}{k(m-1)}}$.  For $a,b$ in $A$, we define $a \circ b:=\omega\cdot a +b$, where $\cdot$ and $+$ are the multiplication and addition operations in $A$ respectively. This is taken to be a left-associative operation. Sometimes we will omit the $\cdot$ for convenience. The binary operation $\circ$ on $A$  induces an $m$-ary operation on $A^{m}$ defined in the following way, 
\begin{equation}\label{eqn:freedefn}
a_1 \circ a_2 \dots \circ a_m := \omega^{m-1} \cdot a_1 + \omega^{m-2} \cdot a_2 + \dots + \omega\cdot a_{m-1} + a_m.
\end{equation}

It is easy to see by direct calculation that the following two lemmas are true. 
\begin{lemma}
The binary operation $\circ$ on $A$ is $k(m-1)$-associative. 
\end{lemma}
\begin{lemma}
The $m$-ary operation on $A^{m}$ induced by the binary operation $\circ$ on $A$ is $k$-associative.
\end{lemma}
Let $X$ be a non-empty set and let $*:X^{m} \rightarrow X$ be an $m$-ary operation. Take $x_1, x_2, \dots, x_n$ in $X$. Recall that there is a bijection between the set of $m$-ary trees on $n$ leaves and the set of $m$-ary parenthesizations of the expression $x_1 * x_2 * \dots * x_n$, see Proposition \ref{HPBijection}. We will write $p_t = p(x_1*x_2 \dots, x_n)_t$ to be the $m$-ary parenthesization of the expression $x_1*x_2*\dots*x_n$ corresponding to the $m$-ary tree $t$. We denote the evaluation of $p_{t}$ with respect to $\circ$ by $p(u_1 \circ u_2 \circ \dots \circ u_n)_{t}.$ When there is no risk of confusion, we omit the subscript $t$. 

\begin{lemma}{\label{multiplication}}
 Let $p(x_1*x_2* \dots*x_n)_t$ be an $m$-ary parenthesization of $x_1 * x_2 * \dots * x_n$ corresponding to the $m$-ary tree on $n$ leaves $t$. Let $( \delta^{l_1}(t), \delta^{l_2}(t),\dots,\delta^{l_{m}}(t) )$ be the depth of $t$. Then we have that\\
 $$p(u_1 \circ u_2 \circ  \dots \circ u_n)_t=\omega^{ \sum\limits_{i=1}^{m} (m-i)\delta^{l_i}_{1}(t)} \cdot u_1 + \omega^{ \sum\limits_{i=1}^{m} (m-i)\delta^{l_i}_{2}(t)} \cdot u_2 + \dots + \omega^{ \sum\limits_{i=1}^{m} (m-i)\delta^{l_i}_{n}(t)} \cdot u_n.$$
 \begin{proof} Recall that $n$ satisfies the equation $n = m + g(m-1)$ for some integer $g \geq 0$. We prove the result by induction on $g$. When $g=0$ there is only on tree to consider, namely $t = \varepsilon \wedge \varepsilon \wedge \dots \wedge \varepsilon$. 
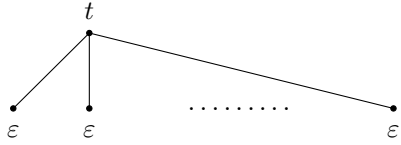
\begin{center}
\begin{tikzpicture}
\draw[fill=black] (0,0) circle (1pt);
\draw[fill=black] (-1,-1) circle (1pt);
\draw[fill=black] (0,-1) circle (1pt);
\draw[fill=black] (4,-1) circle (1pt);

\draw[thin] (0,0) -- (-1,-1);
\draw[thin] (0,0) -- (0,-1);
\draw[thin] (0,0) -- (4,-1);

\node at (0,0.3) {$t$};
\node at (-1,-1.3) {$\varepsilon$};
\node at (0,-1.3) {$\varepsilon$};
\node at (4,-1.3) {$\varepsilon$};
\node at (2,-1) {$\dots \dots \dots $};
\end{tikzpicture} 
 \end{center}
For this tree we have that $\delta^{l_i}_{j} =\delta_{ij}$, where the right hand side is the usual Kronecker delta function. it is easy to see that the statement holds in this case by the definition of $u_1 \circ u_2 \circ \dots \circ u_m$ in (\ref{eqn:freedefn}).

Now suppose that the result holds for $n=m+g^{\prime}(m-1)$ for all $g^{\prime} \leq g$. We consider the $g+1$ case. Let $t$ be an $m$-ary tree with $n= m+(g+1)(m-1)$ leaves. Then we may write $t =  t_1 \wedge \dots \wedge t_m$ where each $t_i$ is the subtree rooted at the $i^{\text{th}}$ child of the root of $t$. Each subtree $t_i$ has $n_i < n$ leaves and $n_1 + n_2 \dots + n_m = n$. In writing $t$ as the meet of its sub-trees at the root, we partition the leaves of $t$. We identify each leaf of $t$ with a tuple $(h,j)$ if it lies in the subtree $t_h$ and it is the $j^{\text{th}}$ leaf in the linear order on the leaves of $t_h$, where $1 \leq j \leq n_{h}$. Therefore for the leaf $(h,j)$ we have that, 
\begin{equation}
\delta^{l_i}_{h,j}(t) =
\begin{cases}
\delta^{l_i}_{h,j}(t_h) +1  & \text{ if } i=h  ;\\
\delta^{l_i}_{h,j}(t_h) & \text{otherwise}.\\
\end{cases}
\end{equation}
From the above equation, it follows that, 
\begin{equation}\label{eqn:proofeqn}
(m-i)\delta^{l_i}_{h,j}(t) =
\begin{cases}
(m-i)\delta^{l_i}_{h,j}(t_h) +(m-i)  & \text{ if } i=h  ;\\
(m-i)\delta^{l_i}_{h,j}(t_h) & \text{otherwise}.\\
\end{cases}
\end{equation}
The identification of the leaves with the tuples $(h,j)$ gives another labelling of the variables $u_s$, where $1 \leq s \leq n$. Since the variable $u_s$ corresponds to the $s^{\text{th}}$ leaf of $t$, and the $s^{\text{th}}$ leaf is identified with $(h,j)$, then we write $u_{(h,j)}$ for $u_s$. So we have that $$p(u_1 \circ u_2 \circ \dots \circ u_n)_t = p(u_{(1,1)} \circ u_{(1,2)} \dots u_{(m,n_m)})_t.$$
It is then easy to see that, 
\begin{align*}
p(u_{(1,1)} \circ u_{(1,2)} \circ \dots \circ u_{(m,n_m)})_t =& p(u_{(1,1)} \circ \dots \circ u_{(1,n_1)})_{t_{1}} \circ p(u_{(2,1)} \circ \dots \circ \dots \circ u_{(2,n_2)})_{t_{2}} \circ \\
& \dots \circ p(u_{(m,1)} \circ \dots \circ u_{(m,n_m)})_{t_{m}}
\end{align*}
$$=\omega^{m-1}p(u_{(1,1)} \circ \dots \circ u_{(1,n_1)})_{t_{1}}+\omega^{m-2}p(u_{(2,1)} \circ \dots \circ \dots \circ u_{(2,n_2)})_{t_{2}} + \dots + p(u_{(m,1)} \circ \dots \circ u_{(m,n_m)})_{t_{m}}.$$
By in the inductive assumption we have that, 
\begin{align*}
p(u_{(h,1)} \circ u_{(h,2)} \circ \dots \circ u_{(h,n_h)})_{t_{h}} =& \omega^{ \sum\limits_{i=1}^{m} (m-i)\delta^{l_i}_{(h,1)}(t_h)} \cdot u_{(h,1)} + \omega^{ \sum\limits_{i=1}^{m} (m-i)\delta^{l_i}_{(h,2)}(t_h)} \cdot u_{(h,2)} + \\
& \dots + \omega^{ \sum\limits_{i=1}^{m} (m-i)\delta^{l_i}_{(h,n_h)}(t_h)} \cdot u_{(h,n_{h})}
\end{align*}
from which it follows that,
\begin{align*}
\omega^{m-h}p(u_{(h,1)} \circ u_{(h,2)} \circ \dots \circ u_{(h,n_h)})_{t_{h}} =\omega^{ \sum\limits_{i=1}^{m} (m-i)\delta^{l_i}_{(h,1)}(t_h) + (m-h)} \cdot u_{(h,1)} + \\ \omega^{ \sum\limits_{i=1}^{m} (m-i)\delta^{l_i}_{(h,2)}(t_h) + (m-h)} \cdot u_{(h,2)} + \dots +\omega^{ \sum\limits_{i=1}^{m} (m-i)\delta^{l_i}_{(h,n_h)}(t_h) + (m-h)} \cdot u_{(h,n_{h})}.
\end{align*}

By equation (\ref{eqn:proofeqn}),
\begin{align*}
\left( \sum\limits_{i=1}^{m} (m-i)\delta^{l_i}_{(h,j)}(t_h) \right)+ (m-h) &= \left( \sum\limits_{\substack{ i = 1 \\i \neq h}}^{m} (m-i)\delta^{l_i}_{(h,j)}(t_h) \right) +(m-h)\delta^{l_{h}}_{(h,j)}(t_h) +(m-h) \\
&= \left( \sum\limits_{\substack{ i = 1 \\i \neq h}}^{m} (m-i)\delta^{l_i}_{(h,j)}(t) \right) +(m-h)\delta^{l_{h}}_{(h,j)}(t)\\
&= \sum\limits_{i=1}^{m} (m-i)\delta^{l_i}_{(h,j)}(t) 
\end{align*}
Therefore, 
\begin{align*}
p(u_{(1,1)} \circ u_{(1,2)} \circ \dots \circ u_{(m,n_m)})_t =& \omega^{ \sum\limits_{i=1}^{m} (m-i)\delta^{l_i}_{(1,1)}(t)} \cdot u_{(1,1)} + \omega^{ \sum\limits_{i=1}^{m} (m-i)\delta^{l_i}_{(1,2)}(t)} \cdot u_{(1,2)} + \\
 &\dots + \omega^{ \sum\limits_{i=1}^{m} (m-i)\delta^{l_i}_{(m,n_m)}(t)} \cdot u_{(m,n_{m})}, 
\end{align*}
as required. This completes the proof.

\end{proof}
\end{lemma}

\begin{theorem}{\label{mainapp}}
Let $p= p(x_1*x_2 \dots, x_n)_t$ and $p^{\prime} = p^{\prime}(x_1*x_2 \dots, x_n)_{t^{\prime}}$ be two $m$-ary parenthesizations of $x_1*x_2*\dots*x_n$ corresponding to the $m$-ary trees on $n$ leaves $t$ and $t^{\prime}$ respectively. Then $p$ and $p'$ are $k$-equivalent with respect to $k$-associativity if and only if, $$p(u_1 \circ u_2 \dots \circ u_n)_{t} = p^{\prime}(u_1 \circ u_2 \dots \circ u_n)_{t^{\prime}}.$$

\begin{proof}
Suppose the parenthesizations $p$ and $p^{\prime}$ are $k$-equivalent. Then $t$ and $t^{\prime}$ are $k$-equivalent also. By Theorem \ref{TreeChar}, \[  \sum\limits_{i=1}^{m-1} (m-i)\delta^{l_{i}}(t) \equiv \sum\limits_{i=1}^{m-1} (m-i)\delta^{l_{i}}(t^{\prime}) \text{ mod } k(m-1). \]
Therefore, $$p(u_1 \circ u_2 \dots \circ u_n)_{t} = p^{\prime}(u_1 \circ u_2 \dots \circ u_n)_{t^{\prime}}$$ by Lemma \ref{multiplication}. 

Suppose that $$p(u_1 \circ u_2 \dots \circ u_n)_{t} = p^{\prime}(u_1 \circ u_2 \dots \circ u_n)_{t^{\prime}},$$ then
$$\omega^{ \sum\limits_{i=1}^{m} (m-i)\delta^{l_i}_{1}(t)} \cdot u_1 + \omega^{ \sum\limits_{i=1}^{m} (m-i)\delta^{l_i}_{2}(t)} \cdot u_2 + \dots + \omega^{ \sum\limits_{i=1}^{m} (m-i)\delta^{l_i}_{n}(t)} \cdot u_n$$ = $$\omega^{ \sum\limits_{i=1}^{m} (m-i)\delta^{l_i}_{1}(t^{\prime})} \cdot u_1 + \omega^{ \sum\limits_{i=1}^{m} (m-i)\delta^{l_i}_{2}(t^{\prime})} \cdot u_2 + \dots + \omega^{ \sum\limits_{r=1}^{m} (m-i)\delta^{l_i}_{n}(t^{\prime})} \cdot u_n.$$ 
Since $u_1,u_2,...,u_n$ are algebraically independent and hence linearly independent in $A$, the coefficients of the $u_i$ on each side of the equation must be equal.\\ Hence,  
$$\omega^{ \sum\limits_{i=1}^{m} (m-i)\delta_{j}^{l_i}(t)} = \omega^{ \sum\limits_{i=1}^{m} (m-i)\delta_{j}^{l_i}(t^{\prime})} \text{ for } 1 \leq j \leq n.$$
 Since $\omega$ has order $k(m-1)$ this implies,\[  \sum\limits_{i=1}^{m-1} (m-i)\delta^{l_{i}}(t) \equiv \sum\limits_{i=1}^{m-1} (m-i)\delta^{l_{i}}(t^{\prime}) \text{ mod } k(m-1). \] Hence $t$ and $t^{\prime}$ are $k$-equivalent by Theorem \ref{TreeChar} which implies that $p$ and $p^{\prime}$ are also $k$-equivalent by Remark \ref{rot&assoc}.
\end{proof}
\end{theorem}
\begin{example}
In example \ref{rotex} we saw that the 3-ary parenthesization $((x_1x_2x_3)x_4x_5)x_6x_7$ is 2-equivalent to $x_1((x_2x_3x_4)x_5x_6)x_7$. Let us check the above theorem for this example. \\
The depth of the first tree is  $$( \delta^{l_1}=(3,2,2,1,1,0,0) , \delta^{l_2}= (0,1,0,1,0,1,0), \delta^{l_3}=(0,0,1,0,1,0,1) ).$$ Therefore the valuation of $((x_1x_2x_3)x_4x_5)x_6x_7$ with respect to $\circ$ is $$\omega^6x_1 + \omega^5x_2 + \omega^4x_3 + \omega^3x_4 + \omega^2x_5 + \omega x_6 +x_7.$$The depth of $x_1((x_2x_3x_4)x_5x_6)x_7$ is $$(\delta^{l_1}=(1,2,1,1,0,0,0) , \delta^{l_2}= (0,1,2,1,2,1,0), \delta^{l_3}=(0,0,0,1,0,1,1)),$$ hence the valuation of $x_1((x_2x_3x_4)x_5x_6)x_7$ with respect to $\circ$ is, 
$$\omega^2x_1 + \omega^5x_2 + \omega^4x_3 + \omega^3x_4 + \omega^2x_5 + \omega x_6 +x_7.$$
Since $\omega$ has order 4 the valuations are equal.
\end{example}

\section{Modular Fuss-Catalan Number}
Recall that we define the $(k)$-modular Fuss-Catalan number $C^{m}_{k,n}$ to be the number of $k$-equivalence classes of parenthesizations of $x_0*x_1*\dots*x_n$. In the previous sections we saw that $k$-associativity corresponds to $k$-rotation and $k$-compression. Therefore $C^{m}_{k,n}$ also counts the $k$-equivalence classes of $(m-1)$-Dyck paths of length $n$. In this section we follow the strategy of \cite[\S 5]{HeinHuangv2} to derive a closed formula for $C^{m}_{k,n}$, see Theorem \ref{closedFormula}. By Proposition \ref{minimalunique}, each $k$-equivalence class has a unique minimal element. Therefore to count the number of $k$-equivalence classes, we just need to count the number of minimal elements. For the rest of this section, we fix integers $m \geq 2$, $g \geq 0$, $k \geq 1$ and $n = m+g(m-1)$.

Assume that $l$ is a positive integer in $\{1,2, \dots , n\}$ such that $(m-1)$ divides $l$. Let $N$ denote the up-step $(1,1)$ and $S$ denote the down-step $(1,-1)$ in $\mathbb{Z}^2$. Denote by $\mathscr{C}^{\prime m}_{k,n,l}$ the set of all strings (lattice paths) of the form $N^{l}SN^{i_1}SU^{i_2} \dots SN^{i_n}$ such that $i_1 + i_2 + \dots + i_n = n-l$ where $(m-1)|i_p$ for all $1 \leq p \leq n$ and  $ 0 \leq i_1 , \dots , i_n <k(m-1)$. So $\mathscr{C}^{\prime m}_{k,n,l}$ is a set of lattice paths of length $n$ where the first up-step is of size $l$. For integers $ 1 \leq j \leq k$, denote by $m_j$ the number of $(j-1)(m-1)$s appearing among the $i_1, i_2, \dots, i_n \in \{0, m-1 , 2(m-1), \dots , (k-1)(m-1) \}$. 

It is easy to see that $$|\mathscr{C}^{\prime m}_{k,n,l}| = \sum_{\substack{m_1 + \dots m_k = n \\ m_2 + 2m_3 + \dots (k-1)m_k = \frac{n-l}{m-1}}} \binom{n}{m_1 , m_2, \dots , m_k}.$$

For a string $w = N^{l}SN^{i_1}SU^{i_2} \dots SN^{i_n}$ in $\mathscr{C}^{\prime m}_{k,n,l}$ and $j$ in $\{0,1, \dots , n-1\}$ we define
$$w^{\bullet j} \vcentcolon= N^{l}SN^{i_{j+1}}S \dots SN^{i_n}SU^{i_1}S \dots SU^{i_j}.$$ 
Let $\mathscr{C}^{m}_{k,n,l}$ be the subset of strings in $\mathscr{C}^{\prime m}_{k,n,l}$ which are $(m-1)$-Dyck paths of length $n$. Since $(m-1)$-Dyck paths can be thought of as $1$-Dyck paths where up-steps come in multiples of $m-1$, the following lemmas follow by similar arguments to Lemma 5.5 and Lemma 5.6 from \cite{HeinHuangv2}, which may be thought of as the $m=2$ case. So we will state the lemmas without proof. 

\begin{lemma}
For  a string $w$ in $\mathscr{C}^{\prime m}_{k,n,l}$ the set $\{ 0 \leq j \leq n-1 : w^{\bullet j} \in \mathscr{C}^{m}_{k,l,n} \}$ has cardinality $l$.
\end{lemma}

Let $\phi$ be the following map, $$\phi : \mathscr{C}^{m}_{k,n,l} \times \{0,1, \dots , n-1 \} \rightarrow \mathscr{C}^{\prime m}_{k,n,l},$$ 
$$ (w,j) \mapsto w^{\bullet j}.$$
 
\begin{lemma}{\label{fiber}}
For a string $w$ in $\mathscr{C}^{\prime m}_{k,n,l}$, the fibre $\phi^{-1}(w)$ of $\phi$ over $w$ has cardinality $|\phi^{-1}(w)| = l$.
\end{lemma}
By Proposition \ref{minimalunique}, the $(k)$-modular Fuss-Catalan number counts the number of minimal $(m-1)$-Dyck paths. Moreover by Proposition \ref{minimal}, minimal Dyck paths $D$ satisfy $d(D)=(d_1, d_2 , \dots ,d_n)$ where $d_i < k(m-1)$ for $i\neq1$. Combining the results of Proposition \ref{minimal}, Proposition \ref{minimalunique} and Lemma \ref{fiber} we have that $$|\mathscr{C}^{m}_{k,n,l}| = \frac{l}{•n}|\mathscr{C}^{\prime m}_{k,n,l}|.$$
Let $\mathscr{C}^{m}_{k,n}$ be the set of minimal $(m-1)-$Dyck paths, then we have that  $$|\mathscr{C}^{m}_{k,n}| = \sum_{\substack{ 1 \leq l \leq n \\ (m-1) | l}} |\mathscr{C}^{m}_{k,n,l}|.$$
Therefore, $$C_{k,n}^{m} = |\mathscr{C}^{m}_{k,n}| = \sum_{\substack{ 1 \leq l \leq n \\ m-1 | l}} { \frac{l}{n} \sum_{\substack{ m_1 + \dots +m_k = n \\ m_2 + 2m_3 + \dots (k-1)m_k = \frac{n-l}{m-1}}}}{\binom{n}{m_1, \dots m_k}}$$
is the number of minimal $(m-1)$-Dyck paths of length $n$, so by Proposition \ref{treePathBij}, it is the number of minimal $m$-ary trees of length $n+1$. This completes the proof of Theorem \ref{closedFormula}. 

\begin{example} In this example, we will count the number of $2$-equivalence classes of $3$-ary trees with $7$ leaves. 
There are twelve 3-ary trees altogether. See the figure below for the complete list.
\begin{figure}[H]
\begin{center}
\begin{tikzpicture}
\draw[fill=black] (-14,0) circle (1pt);
\draw[fill=black] (-13,1) circle (1pt);
\draw[fill=black] (-13,0) circle (1pt);
\draw[fill=black] (-12,0) circle (1pt);
\draw[fill=black] (-12,2) circle (1pt);
\draw[fill=black] (-12,1) circle (1pt);
\draw[fill=black] (-11,1) circle (1pt);
\draw[fill=black] (-11,3) circle (1pt);
\draw[fill=black] (-11,2) circle (1pt);
\draw[fill=black] (-10,2) circle (1pt);

\draw[thin] (-14,0) -- (-13,1);
\draw[thin] (-13,1) -- (-13,0);
\draw[thin] (-13,1) -- (-12,0);
\draw[thin] (-13,1) -- (-12,2);
\draw[thin] (-12,2) -- (-12,1);
\draw[thin](-12,2) -- (-11,1);
\draw[thin] (-12,2) -- (-11,3);
\draw[thin] (-11,3) -- (-11,2);
\draw[thin] (-11,3) -- (-10,2);

\node at (-11,3.25) {$T_{1}$};

\draw[fill=black] (-9,0) circle (1pt);
\draw[fill=black] (-8,0) circle (1pt);
\draw[fill=black] (-7,0) circle (1pt);
\draw[fill=black] (-8,1) circle (1pt);
\draw[fill=black] (-7,1) circle (1pt);
\draw[fill=black] (-6,1) circle (1pt);
\draw[fill=black] (-7,2) circle (1pt);
\draw[fill=black] (-6,2) circle (1pt);
\draw[fill=black] (-8,2) circle (1pt);
\draw[fill=black] (-7,3) circle (1pt);

\draw[thin] (-8,1) -- (-8,0);
\draw[thin] (-8,1) -- (-9,0);
\draw[thin] (-8,1) -- (-7,0);
\draw[thin] (-7,2) -- (-8,1);
\draw[thin] (-7,2) -- (-7,1);
\draw[thin] (-7,2) -- (-6,1);
\draw[thin] (-7,3) -- (-8,2);
\draw[thin] (-7,3) -- (-7,2);
\draw[thin] (-7,3) -- (-6,2);

\node at (-7,3.25) {$T_{2}$};

\draw[fill=black] (-4,1) circle (1pt);
\draw[fill=black] (-5,0) circle (1pt);
\draw[fill=black] (-4,0) circle (1pt);
\draw[fill=black] (-3,0) circle (1pt);
\draw[fill=black] (-3,2) circle (1pt);
\draw[fill=black] (-3,1) circle (1pt);
\draw[fill=black] (-2,1) circle (1pt);
\draw[fill=black] (-4,2) circle (1pt);
\draw[fill=black] (-5,2) circle (1pt);
\draw[fill=black] (-4,3) circle (1pt);

\draw[thin] (-4,3) -- (-5,2);
\draw[thin] (-4,3) -- (-4,2);
\draw[thin] (-4,3) -- (-3,2);
\draw[thin] (-3,2) -- (-3,1);
\draw[thin] (-3,2) -- (-2,1);
\draw[thin] (-3,2) -- (-4,1);
\draw[thin] (-4,1) -- (-5,0);
\draw[thin] (-4,1) -- (-4,0);
\draw[thin] (-4,1) -- (-3,0);

\node at (-4,3.25) {$T_{3}$};

\end{tikzpicture}
\end{center}

\begin{tikzpicture}
\draw[fill=black] (0,0) circle (1pt);
\draw[fill=black] (1,0) circle (1pt);
\draw[fill=black] (2,0) circle (1pt);
\draw[fill=black] (1,1) circle (1pt);
\draw[fill=black] (0,1) circle (1pt);
\draw[fill=black] (2,1) circle (1pt);
\draw[fill=black] (1,2) circle (1pt);
\draw[fill=black] (2,2) circle (1pt);
\draw[fill=black] (3,2) circle (1pt);
\draw[fill=black] (2,3) circle (1pt);

\draw[thin] (1,1) -- (0,0);
\draw[thin] (1,1) -- (2,0);
\draw[thin] (1,1) -- (1,0);
\draw[thin] (1,2) -- (2,1);
\draw[thin] (1,2) -- (1,1);
\draw[thin] (1,2) -- (0,1);
\draw[thin] (2,3) -- (1,2);
\draw[thin] (2,3) -- (2,2);
\draw[thin] (2,3) -- (3,2);

\node at (2,3.25) {$T_{4}$};

\draw[fill=black] (4,0) circle (1pt);
\draw[fill=black] (5,0) circle (1pt);
\draw[fill=black] (6,0) circle (1pt);
\draw[fill=black] (5,1) circle (1pt);
\draw[fill=black] (4,1) circle (1pt);
\draw[fill=black] (6,1) circle (1pt);
\draw[fill=black] (4,2) circle (1pt);
\draw[fill=black] (5,2) circle (1pt);
\draw[fill=black] (6,2) circle (1pt);
\draw[fill=black] (5,3) circle (1pt);

\draw[thin] (5,1) -- (4,0);
\draw[thin] (5,1) -- (5,0);
\draw[thin] (5,1) -- (6,0);
\draw[thin] (5,2) -- (4,1);
\draw[thin] (5,2) -- (5,1);
\draw[thin] (5,2) -- (6,1);
\draw[thin] (5,3) -- (4,2);
\draw[thin] (5,3) -- (5,2);
\draw[thin] (5,3) -- (6,2);

\node at (5,3.25) {$T_{5}$};

\draw[fill=black] (8,0) circle (1pt);
\draw[fill=black] (9,0) circle (1pt);
\draw[fill=black] (10,0) circle (1pt);
\draw[fill=black] (9,1) circle (1pt);
\draw[fill=black] (8,1) circle (1pt);
\draw[fill=black] (10,1) circle (1pt);
\draw[fill=black] (9,2) circle (1pt);
\draw[fill=black] (8,2) circle (1pt);
\draw[fill=black] (7,2) circle (1pt);
\draw[fill=black] (9,2) circle (1pt);
\draw[fill=black] (8,3) circle (1pt);

\draw[thin] (9,1) -- (8,0);
\draw[thin] (9,1) -- (9,0);
\draw[thin] (9,1) -- (10,0);
\draw[thin] (9,2) -- (8,1);
\draw[thin] (9,2) -- (9,1);
\draw[thin] (9,2) -- (10,1);
\draw[thin] (8,3) -- (7,2);
\draw[thin] (8,3) -- (8,2);
\draw[thin] (8,3) -- (9,2);

\node at (8,3.25) {$T_{6}$};

\draw[fill=black] (12,0) circle (1pt);
\draw[fill=black] (13,0) circle (1pt);
\draw[fill=black] (14,0) circle (1pt);
\draw[fill=black] (13,1) circle (1pt);
\draw[fill=black] (12,1) circle (1pt);
\draw[fill=black] (11,1) circle (1pt);
\draw[fill=black] (12,2) circle (1pt);
\draw[fill=black] (13,2) circle (1pt);
\draw[fill=black] (14,2) circle (1pt);
\draw[fill=black] (13,3) circle (1pt);

\draw[thin] (13,1) -- (12,0);
\draw[thin] (13,1) -- (13,0);
\draw[thin] (13,1) -- (14,0);
\draw[thin] (12,2) -- (13,1);
\draw[thin] (12,2) -- (12,1);
\draw[thin] (12,2) -- (11,1);
\draw[thin] (13,3) -- (12,2);
\draw[thin] (13,3) -- (13,2);
\draw[thin] (13,3) -- (14,2);

\node at (13,3.25) {$T_{7}$};
\end{tikzpicture}

\begin{center}
\begin{tikzpicture}
\draw[fill=black] (0,0) circle (1pt);
\draw[fill=black] (1,0) circle (1pt);
\draw[fill=black] (2,0) circle (1pt);
\draw[fill=black] (1,1) circle (1pt);
\draw[fill=black] (0,1) circle (1pt);
\draw[fill=black] (-1,1) circle (1pt);
\draw[fill=black] (0,2) circle (1pt);
\draw[fill=black] (1,2) circle (1pt);
\draw[fill=black] (-1,2) circle (1pt);
\draw[fill=black] (0,3) circle (1pt);

\draw[thin] (1,1) -- (1,0);
\draw[thin] (1,1) -- (2,0);
\draw[thin] (1,1) -- (0,0);
\draw[thin] (0,2) -- (0,1);
\draw[thin] (0,2) -- (1,1);
\draw[thin] (0,2) -- (-1,1);
\draw[thin] (0,3) -- (0,2);
\draw[thin] (0,3) -- (1,2);
\draw[thin] (0,3) -- (-1,2);

\node at (0,3.25) {$T_{8}$};

\draw[fill=black] (5,0) circle (1pt);
\draw[fill=black] (6,0) circle (1pt);
\draw[fill=black] (7,0) circle (1pt);
\draw[fill=black] (6,1) circle (1pt);
\draw[fill=black] (5,1) circle (1pt);
\draw[fill=black] (4,1) circle (1pt);
\draw[fill=black] (5,2) circle (1pt);
\draw[fill=black] (4,2) circle (1pt);
\draw[fill=black] (3,2) circle (1pt);
\draw[fill=black] (4,3) circle (1pt);

\draw[thin] (6,1) -- (5,0);
\draw[thin] (6,1) -- (6,0);
\draw[thin] (6,1) -- (7,0);
\draw[thin] (5,2) -- (6,1);
\draw[thin] (5,2) -- (5,1);
\draw[thin] (5,2) -- (4,1);
\draw[thin] (4,3) -- (4,2);
\draw[thin] (4,3) -- (5,2);
\draw[thin] (4,3) -- (3,2);

\node at (4,3.25) {$T_{9}$};

\draw[fill=black] (8,0) circle (1pt);
\draw[fill=black] (9,0) circle (1pt);
\draw[fill=black] (10,0) circle (1pt);
\draw[fill=black] (11,0) circle (1pt);
\draw[fill=black] (12,0) circle (1pt);
\draw[fill=black] (13,0) circle (1pt);
\draw[fill=black] (9,1) circle (1pt);
\draw[fill=black] (12,1) circle (1pt);
\draw[fill=black] (13,1) circle (1pt);
\draw[fill=black] (11,3) circle (1pt);

\draw[thin] (9,1) -- (8,0);
\draw[thin] (9,1) -- (9,0);
\draw[thin] (9,1) -- (10,0);
\draw[thin] (12,1) -- (11,0);
\draw[thin] (12,1) -- (12,0);
\draw[thin] (12,1) -- (13,0);
\draw[thin] (11,3) -- (9,1);
\draw[thin] (11,3) -- (12,1);
\draw[thin] (11,3) -- (13,1);

\node at (11,3.25) {$T_{10}$};
\end{tikzpicture}
\end{center}

\begin{center}
\begin{tikzpicture}
\draw[fill=black] (0,0) circle (1pt);
\draw[fill=black] (1,0) circle (1pt);
\draw[fill=black] (2,0) circle (1pt);
\draw[fill=black] (3,0) circle (1pt);
\draw[fill=black] (4,0) circle (1pt);
\draw[fill=black] (5,0) circle (1pt);
\draw[fill=black] (1,1) circle (1pt);
\draw[fill=black] (2.5,1) circle (1pt);
\draw[fill=black] (4,1) circle (1pt);
\draw[fill=black] (2.5,3) circle (1pt);

\draw[thin] (1,1) -- (0,0);
\draw[thin] (1,1) -- (1,0);
\draw[thin] (1,1) -- (2,0);
\draw[thin] (4,1) -- (3,0);
\draw[thin] (4,1) -- (4,0);
\draw[thin] (4,1) -- (5,0);
\draw[thin] (2.5,3) -- (1,1);
\draw[thin] (2.5,3) -- (2.5,1);
\draw[thin] (2.5,3) -- (4,1);

\node at (2.5,3.25) {$T_{11}$};

\draw[fill=black] (7,0) circle (1pt);
\draw[fill=black] (8,0) circle (1pt);
\draw[fill=black] (9,0) circle (1pt);
\draw[fill=black] (10,0) circle (1pt);
\draw[fill=black] (11,0) circle (1pt);
\draw[fill=black] (12,0) circle (1pt);
\draw[fill=black] (8,1) circle (1pt);
\draw[fill=black] (7,1) circle (1pt);
\draw[fill=black] (11,1) circle (1pt);
\draw[fill=black] (9.5,3) circle (1pt);

\draw[thin] (8,1) -- (7,0);
\draw[thin] (8,1) -- (8,0);
\draw[thin] (8,1) -- (9,0);
\draw[thin] (11,1) -- (10,0);
\draw[thin] (11,1) -- (11,0);
\draw[thin] (11,1) -- (12,0);
\draw[thin] (9.5,3) -- (8,1);
\draw[thin] (9.5,3) -- (11,1);
\draw[thin] (9.5,3) -- (7,1);

\node at (9.5,3.25) {$T_{12}$};

\end{tikzpicture}
\end{center}

  \caption{The complete list of the 3-ary trees with 7 leaves.}
  \label{fig:3Trees7}
\end{figure}
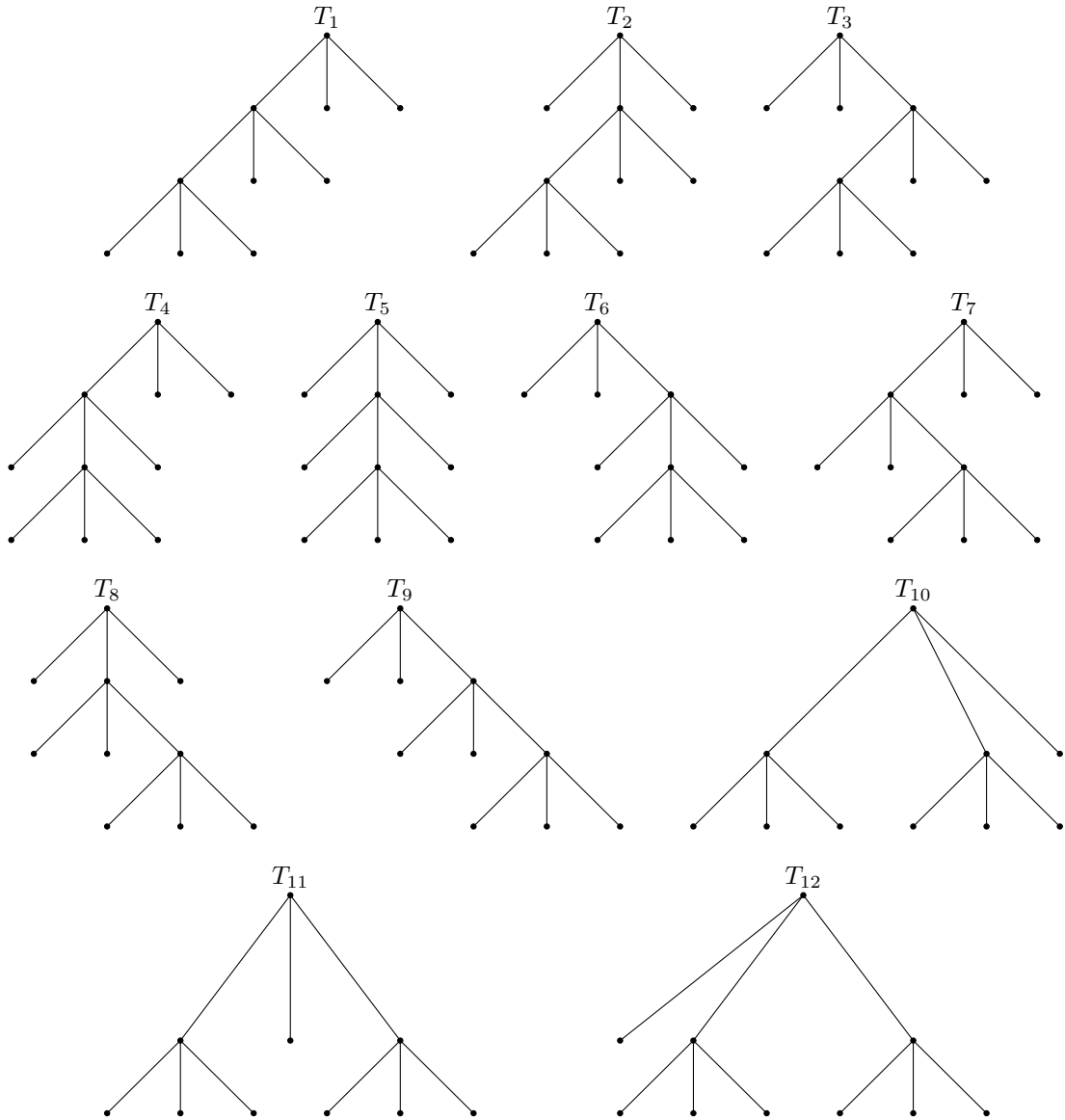

Observe that the trees $T_{1}, T_{2}$ and $T_{3}$ in the top row correspond to the following parenthesizations  $$((x_1*x_2*x_3)*x_4*x_5)*x_6*x_7, $$ $$x_1*((x_2*x_3*x_4)*x_5*x_6)*x_7,$$ 
$$x_1*x_2*((x_3*x_4*x_5)*x_6*x_7)$$
respectively. We can see that we get the tree $T_{2}$ from the tree $T_{1}$ by a 2-rotation at the root of $T_{1}$, and likewise we get the tree $T_{3}$ from the tree $T_{2}$ by a 2-rotation at the root of $T_{2}$. Therefore $T_{1}, T_{2}$ and $T_{3}$ belong to the same $2$-equivalence class. Further observe that the other trees cannot be $2$-rotated because they don't contain a subtree of form required to perform a $2$-rotation. So we conclude that $C_{2,6}^{3}=10$. Let us check this against the closed formula of Theorem \ref{closedFormula}. 
\begin{align*}
C^{3}_{2,6} &= \sum_{\substack{ 1 \leq l \leq 6 \\2 | l}} { \frac{l}{6} \sum_{\substack{ m_1 + m_2= 6\\ m_2= \frac{6-l}{2}}}}{\binom{6}{m_1,m_2}} \\
&= { \frac{2}{6}\binom{6}{4,2}} + { \frac{4}{6}\binom{6}{5,1}} + { \frac{6}{6}\binom{6}{6,0}} \\
&=  \frac{1}{3}(15) + \frac{2}{3}(6) + 1(1) \\
&=  10.
\end{align*}

\end{example}
 
\bibliographystyle{plain}
 \bibliography{ModularFussCatalanNumbers}
 
\end{document}